\documentclass[journal]{IEEEtran}
\usepackage{subfigure}
\usepackage{multirow}

\usepackage{hyperref}

\usepackage{color}

\usepackage{cite}
\usepackage{pdfsync}
\usepackage{graphicx,epstopdf}

\usepackage{color}

\usepackage{empheq}

\newcommand{\qed}{\hfill \ensuremath{\Box}}

\usepackage{algorithmic}

\usepackage{float}
\floatstyle{ruled}
\newfloat{algorithm}{htbp}{loa}
\floatname{algorithm}{Algorithm}

\usepackage{multirow}
\usepackage{cite}
\usepackage{amsmath, amsthm, amssymb,graphicx,bbm}
\usepackage{psfrag}
\usepackage{bm}

\newcommand{\R}{\mathbb{R}}

\newcommand{\be}{\beta}

\newcommand{\bd}{\mathbf}

\newcommand{\tl}{\tilde}

\newcommand{\x}{\bd{x}}

\newcommand{\ov}{\overline}
\newcommand{\ul}{\underline}

\newcommand{\C}{\mathbb{C}}

\newcommand{\sdp}{\succcurlyeq}
\newcommand{\bma}{\begin{bmatrix}}
\newcommand{\ebma}{\end{bmatrix}}

\newcommand{\nn}{\nonumber}
\newcommand{\te}{\theta}

\newcommand{\mc}{\mathcal}

\DeclareMathOperator{\rank}{rank}
\DeclareMathOperator{\Tr}{Tr}

\DeclareMathOperator{\diag}{diag}
\DeclareMathOperator{\Real}{Re}

\DeclareMathOperator{\convhull}{convhull}

\newtheorem{thm}{Theorem}

\newtheorem{lem}[thm]{Lemma}

\begin{document}
\title{An Optimal and Distributed Method for  Voltage Regulation  in Power Distribution Systems}
\author{Baosen~Zhang,~\IEEEmembership{Member,~IEEE,}
Albert~Y.S.~Lam,~\IEEEmembership{Member,~IEEE,}        \\
        Alejandro~Dom\'{i}nguez-Garc\'{i}a,~\IEEEmembership{Member,~IEEE,}
        and David~Tse,~\IEEEmembership{Fellow,~IEEE}
\thanks{B. Zhang is with the Department of Civil and Environmental Engineering and the Department of Management and Science Engineering at Stanford University; E-mail: zhangbao@stanford.edu. A.Y.S. Lam is with the Department of Computer Science at Hong Kong Baptist University; E-mail: albertlam@ieee.org. A. Dom\'{i}nguez-Garc\'{i}a is with the   Department of Electrical and Computer Engineering of the University of Illinois at Urbana-Champaign;  E-mail: aledan@illinois.edu. D. Tse is  with the  Department of Electrical Engineering  at Stanford University; E-mail:  dntse@stanford.edu}.
}%
 \vspace{-0.1in}
 \maketitle
 
 \vspace{-0.1in}

\begin{abstract}
This paper addresses the problem of voltage regulation in power distribution networks with deep-penetration of distributed energy resources, e.g., renewable-based generation, and storage-capable loads such as plug-in hybrid electric vehicles. We cast the problem as an optimization program, where the objective is to minimize the losses in the network subject to constraints on bus voltage magnitudes,  limits on active and reactive power injections,  transmission line thermal limits and losses. We provide sufficient conditions under which the optimization problem can be solved via its convex relaxation. Using data from existing networks, we show that these sufficient conditions are expected to be satisfied by most networks. We also provide an efficient distributed algorithm to solve the problem.  { The algorithm adheres to a communication topology  described by a graph  that is the same as the graph that describes the electrical network topology. We illustrate the operation of the algorithm, including its robustness against communication link failures, through several case studies involving 5-,   34- and 123-bus  power distribution systems.}

\end{abstract}


\section{Introduction}
Electric power distribution systems will undergo radical transformations in structure and functionality due to the advent of initiatives like the  US DOE \textit{Smart Grid} \cite{smart_grid}, and its European counterpart \textit{Electricity Networks of the Future} \cite{smart_grid_eu}. These transformations are enabled by the integration of (i) advanced communication and control,  (ii)  renewable-based variable generation resources, e.g.,  photovoltaics (PVs), and  (iii) new storage-capable loads, e.g.,  plug-in hybrid electric vehicles (PHEVs). These distributed generation and storage resources are commonly referred to as distributed energy resources (DERs).
It has been acknowledged (see, e.g., \cite{Carvalho:2008}) that massive penetration of DERs in distribution networks is likely to cause voltage regulation problems due to the fact that typical values of transmission line resistance-to-reactance, $r/x$, ratios are such that bus voltage magnitudes are fairly sensitive to variations in active power injections (see, e.g.,  \cite{Keane:2011}).
Similarly massive penetration of PHEVs can potentially create substantial voltage drops ~\cite{GuGr:09}. The objective of this paper is to address the problem of voltage regulation in electric power distribution networks with deep penetration of DERs.

{As of today, voltage regulation in distribution networks is accomplished through  tap-changing under-load transformers, set voltage regulators, and fixed/switched capacitors. While these devices---the operation of which is mechanical in nature---are effective in managing slow variations (on the time-scale of hours) in voltage, their lifetime could be dramatically reduced from the increased number of operations needed to handle faster voltage variations due to sudden changes (on the time-scale of minutes) in active power generated or consumed  by DERs. An alternative to the use of these   voltage regulation devices for handling fast variations is to utilize the power electronics interfaces of the DERs themselves.}

While active power control is the primary function of these  interfaces, when being properly controlled, they can also provide reactive power. Thus, they provide a mechanism to control reactive power injections, which in turn can be used for voltage control  (see, e.g., \cite{JoOo:00,LoKr:01}).   In this regard, there are  existing PV rooftop and pole-mount   solutions  that  provide such functionality (see, e.g., \cite{petra,solar_bridge}). These solutions  are endowed with wireless \cite{petra}, and power-line communications \cite{solar_bridge}, which is   a key for  controlling a large number of devices without overlaying a separate communication network. Additionally, as noted earlier, bus voltages in a distribution network are sensitive to changes in active power injections. Thus, storage-capable DERs,  and demand response resources (DRRs), provide a second voltage control mechanism as they can be used, to some extent, to shape active power injections.


This paper  proposes a method for voltage regulation   in distribution networks that  relies on the utilization of reactive-power-capable DERs, and to some extent, on the control of active power injections enabled by storage-capable DERs and DRRs. {This method is intended to supplement the action of conventional voltage regulation devices,  while i) minimizing their usage by handling faster voltage variations due to changes in renewable-based power injections, and  ii) having them intervene only during extreme circumstances rather than minor, possibly temporary,  violations. In this regard, in subsequent developments, we assume that there is a separation in the (slow) time-scale  in which the settings of   conventional voltage regulation devices are adjusted and the (fast) time-scale in which our proposed method operates. 
For instance, the settings of conventional devices can   be optimized every hour in anticipation of the overall change in load (e.g., air conditioners being turned on in a hot afternoon). Then, within each hour, our  method is  utilized  to regulate voltage  in response to fast variations in  DER-based power injections. A more detailed description of the ideas above is provided in  Section~\ref{voltage_regulation_problem}.} 

{The voltage regulation problem can be cast as a optimization program where the objective is to minimize network losses\footnote{{Any objective function that is strictly increasing can be used and all of the results in this paper remains unchanged. We focus on loss minimization due to its relevance for long-run economic savings.}}  subject to (i) constraints on bus voltage magnitudes, (ii) upper and lower limits on active and reactive power injections, (iii) upper limits on transmission line flows, and (iv) upper limits on transmission line losses. The decision variables are the bus voltages; however the actual control mechanism to fix these are  reactive (and to some extent active) power injections in  the  bus of the network.}

{The contributions of this paper are two fold. First, we establish sufficient conditions under which the voltage regulation problem can be solved via an equivalent semidefinite programming (SDP) problem, thus convexifying the original problem. This equivalence also leads to a simple method to establish whether or not the original problem has a feasible solution based on the solution of the convexified problem.  However, it is important to note that even if the convexified SDP  provides the solution to the original optimization problem, existing  algorithms for solving SDPs  are not computationally efficient for solving large problems \cite{Boyd04}. Therefore, these algorithms are not practical for realizing our ideas in a realistic power distribution network with thousands of buses. Furthermore, even if there is a centralized solver with sufficient computing power, the communication infrastructure may not be able to reliably transmit all problem data to a centralized  location with small enough delays. This is where the second contribution of our work lies---the  development of a distributed algorithm   for efficiently  solving convexified SDPs on tree networks; as it will shown  later, a key feature of this algorithm is that is robust against communication   failures.}

Previous works  that  addressed the voltage regulation problem in distribution networks also cast  it as an optimization problem; however, in contrast with our work, the solution methods proposed in these  earlier works are sub-optimal, and in most cases, they rely on a centralized decision maker that has access to all the data defining the problem \cite{Horak10,Rumley08}. For example   in \cite{Ba:07}, the authors  propose  the use of reactive-power-capable DERs for voltage control and  the objective is to minimize the DER reactive power contributions subject to the power flow equations and other constraints. However, the solution approach proposed in  \cite{Ba:07}  relies on linearizing the power flow equations around some operating point, rendering a linear program; therefore,  this approach provides a  sub-optimal   solution. In the same vein, and although the objective function is different to the one considered in \cite{Ba:07}, the solution proposed in \cite{Turitsyn:2010}  also relies on linearization.

Other optimization-based  approaches to address the voltage regulation problem in distribution networks rely on optimal power flow (OPF) solvers   developed for transmission networks (see, e.g., \cite{Villacci:2006,AL11}). For instance, in \cite{AL11}, the authors use a Newton-type method to solve the Lagrange dual of the optimization program they consider in their problem; however, since the primal problem is, in general, not convex, there is no guarantee that the solution of the dual problem is  globally optimal, or even  physically meaningful.

{The key  to solve the the voltage regulation problem is to wite it as a rank-constrained SDP \cite{Bai08}, where the decision variable is a positive semidefinite matrix constrained to have rank 1, which, in general, makes the problem not convex. The problem can be convexified by dropping the rank-1 constraint---the conundrum is then to establish when the solution of the convexified problem also provides a global solution to the original non-convex problem. 

Recently there has been a sequence of papers on attempting to answer the question above spurred by the observation in \cite{Lavaei12} that the convex relaxation observed above is tight for many IEEE benchmark transmission networks. Several independent works \cite{Zhang11a,Lavaei11c,Bose11} provided a partial answer: the convex relaxation is tight if the network has a {\em tree} topology and certain constraints on the bus power injections are satisfied. All these results, which are particularly relevant for distribution networks, were unified and strengthened in \cite{LTZ12} through a investigation of the underlying geometry of the optimization problem.} It is important to note that,  even in the situations where the convex relaxation is tight, there might be multiple local optimal solutions and local search algorithms may not converge to the globally optimal one (see \cite{Lavaei11c} for a further discussion on this). Finally, it is well known (see, e..g, \cite{Stott09}) that solution methods to OPF-type problems for transmission networks   tend   to perform poorly  in  distribution networks due to the $r/x$ values.

An independent---but related---work has recently appeared in   \cite{Farivar12}. Although, a direct comparison of both works is difficult due to the difference in assumptions and constraints (for example \cite{Farivar12} does not require all buses to have DERs), the solution method proposed by the authors in  \cite{Farivar12} is also globally optimal; however, their  solution method requires a centralized processor, whereas as it is shown later,  our solution method is  amenable for a distributed implementation.  With respect to this, recent independent works to ours are the ones in \cite{Kraning12} and \cite{Li12b}; where the authors in \cite{Kraning12} propose a distributed algorithm to solve convex relaxations for OPF-type problems in general networks and the authors in \cite{Li12b} considered demand response in the distribution network.  

This paper builds on the results of \cite{Zhang11a,LTZ12} and extends them by taking into account the reactive power injections, and considering tight voltage magnitude constraints. Previous results either: (i) ignore reactive power (e.g. \cite{Zhang11,LTZ12,Gan12,Li12} ), (ii) assume there are not active lower bounds (e.g. \cite{Lavaei12,Li12b}), or (iii) assume that there are no voltage upper bounds (e.g \cite{Gan12}).



The remainder of this paper is organized as follows.  Section~\ref{formulation} discusses the voltage regulation problem and formulates its solution as   an optimization problem. Section \ref{sec:main} states the main theoretical result of the paper, while  Section \ref{sec:geometry} provides a sketch of the proof (the complete proof is provided in Appendix). 
Section~\ref{sec:algo} proposes a distributed algorithm to solve the optimization problem, the performance of which is  illustrated in Section~\ref{sec:simulation} via  case studies. Concluding remarks are presented in Section~\ref{sec:con}.

\section{Preliminaries and Problem Formulation} \label{formulation}
This section introduces the model for the class of power distribution systems considered in this work; in the process,   relevant notations used throughout the paper are also introduced. Subsequently, we discuss the problem of voltage regulation in distribution networks with deep DER penetration. Additionally, we articulate a potential solution that coordinates (a) the utilization of  conventional voltage regulation devices for handling slower voltage variations, and (b) the use of  DERs and DRRs for handling faster variations. The section concludes with the formulation of an optimization problem that enables the realization of (b) above.

\subsection{Power System Distribution  Model}
Consider a power distribution network with $n$ buses. As of today, such networks are mostly radial with a single source of power injection referred to as the feeder (see, e.g., \cite{Kersting06}). Thus, the network topology can be described by  a {\em connected tree}, the edge set of which is denoted by $\mathcal{E}$, where $(i,k) \in \mc{E}$ if $i$ is connected to $k$ by a  transmission line; we write $i\sim k$ if bus $i$ is   connected to bus $k$ and $i \not\sim k$ otherwise, and  $(i,k)$ to denote a transmission line connected between buses  $i$ and $k$.  {Typical examples of power distribution networks with such tree topologies are the IEEE 34- and 123-bus  distribution systems; the topologies of these systems, which are used in the case studies of Section~\ref{sec:simulation}, are displayed in Figs.~\ref{fig:feeder34} and \ref{fig:123bus} respectively.}


Let $V_i=|V_i|\angle \theta_i$ denote bus $i$ voltage, and define the corresponding bus voltage vector  $\bd{v}=[V_1 \; V_2 \; \cdots \; V_n]^T \in \mathbb{C}^n$. Similarly, let $P_i$ and $Q_i$ denote, the active and reactive power injections in bus $i$ respectively, and define the corresponding active and reactive power injection vectors   $\bd{p}=[P_1\; P_2 \; \cdots \; P_n]^T$, $\bd{q}=[Q_1\; Q_2 \; \cdots \;Q_n]^T$. Let $y_{ik}=g_{ik}-jb_{ik}$, with  $b_{ik}>g_{ik}>0$, denote the admittance\footnote{We  adopt the standard assumption that, in normal operating conditions, lines are inductive and the inductive effects dominate resistive effects (see, e.g., \cite{BeVi:00}. Additionally, although the  convention is to write $y_{ik}=g_{ik}+jb_{ik}$,  we chose to  flip the sign of the imaginary part  as it simplifies subsequent developments.} of line $(i,k)$, and let $y_{ii}=j b_{ii}$ denote bus $i$ shunt admittance. Then, the  power flow equations can be compactly written   as
\begin{align} \label{eqn:pq}
\bd{p}+j \bd{q}=\text{Re}\{\text{diag}(\bd{v}\bd{v}^H \bd{Y}^H)\}+ j \text{Im}\{\diag(\bd{v}\bd{v}^H \bd{Y}^H)\},
\end{align}
where $\bd{Y}=[Y_{ik}]$ is the bus admittance matrix ( we use $\bd{A}[i,k]$ to denote the $(i,k)^{th}$ entry of a matrix $\bd{A}$),   $\bd{v}^H$ ($\bd{Y}^H$) denotes the  Hermitian transpose of $\bd{v}$ ($\bd{Y}$), and $\text{diag}(\cdot)$ returns the diagonal of a square matrix as a column vector. The active power flow  through each transmission line  $(i,k)$ is given by
\begin{align}
&P_{ik}=|V_i|^2 g_{ik}+ |V_i||V_k|[b_{ik} \sin(\theta_{ik})-g_{ik} \cos(\theta_{ik})],   \label{eqn:Pik}
\end{align}
whereas the reactive power flow through each transmission line  $(i,k)$ is given by
\begin{align}
&Q_{ik}=|V_i|^2 b_{ik} - |V_i||V_k|[g_{ik} \sin(\theta_{ik})+b_{ik} \cos(\theta_{ik})], \label{eqn:Qik}
\end{align}
where $\theta_{ik}:=\theta_i - \theta_k$.

\subsection{Voltage Control in Networks with Deep DER Penetration} \label{voltage_regulation_problem}
{ 
The objective of the paper is to address the problem of    voltage regulation in   power distribution networks with deep penetration of DERs; specifically, the focus is on the problem of mitigating  voltage variability across the network due to  to fast  (and uncontrolled) changes in the active generated or consumed by DERs. To this end, we rely on i) the use of the power electronics interfaces of the DERs to locally provide some limited amount of reactive power; and ii) to some extent, on the use of storage-capable DERs and DRRs to locally provide (or consume) some amount of active power. In other words, we have a limited ability to shape the active/reactive power injection profile. With respect to this, it is important to note that this ability to shape the active/reactive power injection profile, which in turn will allow us to regulate voltage across the network, and it is intended to supplement the action of conventional voltage regulation devices (e.g., tap-changing under-load transformers, set voltage regulators, and fixed/switched capacitors). 

\begin{figure}[t!] 
\centering
\includegraphics[scale=0.9]{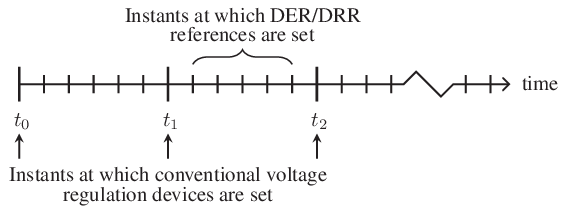}
\caption{Time-scale separation between the  instants in which the settings of conventional voltage regulation devices are decided, and the instants in which the  reference   of  DER and DRRs are set. \label{fig:timescale}}
 \vspace{-0.2in}
\end{figure}

In practice, in order to realize the  ideas above, we envision a hierarchical control architecture; there is a separation in the (slow) time-scale in which the settings of conventional voltage regulation devices are adjusted via the solution to some optimization problem,  and the (fast) time-scale in which voltage regulation through active/reactive power injection shaping is accomplished.  Then, given that fast (and uncontrolled) changes in the DERs active generation (consumption) might cause the voltage to deviate from this reference voltage, a second optimization is performed at regular intervals (e.g., every minute). The timeframe in which the settings of  conventional devices are decided and the reference setting of DERs/DRRs  is graphically depicted in Fig. \ref{fig:timescale}. The solution of this minute-by-minute optimization will provide the  amount of active/reactive power that needs to be locally produced or consumed so as to track the voltage reference. In order words, the minute-by-minute optimization provides the reference values for the amount of active/reactive power to be collectively provided (or consumed) on each bus of the network within the next minute by reactive-power-capable and/or storage-capable DERs and DRRs. These reference values are then passed to the DERs and DRRs local controllers, which will adjust their output accordingly---note that the time-scale in which DER/DRR local controllers act (on the order of milliseconds (see, e.g., \cite{Zou12,Brea10}),   is much faster than the minute-to-minute optimization. Here, it is important to note that the DERs and DRRs are only attempting to correct voltage deviation from the nominal value due to variations in power injections around the power injection profile used to set the conventional voltage regulation devices.

\subsection{Voltage Regulation via DERs/DRRs: Problem Formulation} \label{optimization_formulation}

As stated earlier, the focus of this paper is on developing mechanisms to mitigate  voltage variability across the network due to   fast  (and uncontrolled) changes in the active generated or consumed by DERs; thus, subsequent developments only deal with the inter-hour minute-by-minute optimization mentioned above. As argued in Section~\ref{voltage_regulation_problem}, we assume that at the beginning of each hour, the settings of conventional voltage regulation devices are optimized, which in turns prescribe the values that  each individual bus voltage can take to some   voltage reference $V_i^{ref}$. In order
to achieve the voltage regulation goal above, we rely on a
limited ability to  locally produce and/or consume some limited amount of active/reactive power, and cast the voltage regulation problem as an
optimization program with the objective of minimizing network losses.
}

%

Let $L_{ik}(V_i,V_k)=P_{ik}+P_{ki}$. The total losses in the network are given by $L(\bd{v}):=\sum_{i,k:(i,k)\in \mathcal{E}}L_{ik}(V_i,V_k)$, and the voltage regulation problem can be formulated as
\begin{subequations}
\label{eqn:prob1}
\begin{align}
\underset{\bd{v}~~}{\text{min~~}}  & L(\bd{v})~~ \\
\text{s.t.~~} & |V_i| ={V}_i^{ref}, \forall i	\label{eqn:vol}    \\
 & \underline{P}_i\leq P_i \leq \overline{P}_i, \forall i \label{eqn:p}	\\
 & \underline{Q}_i\leq Q_i \leq \overline{Q}_i, \forall i \label{eqn:q}	\\
 & |P_{ik}| \leq \ov{P}_{ik}, \; \forall i \sim k \label{eqn:pik}\\
 & L_{ik}(V_i,V_k) \leq \ov{L}_{ik}, \; \forall i \sim k  \label{eqn:lik} \\
 & { P_i = \sum_{k \sim i} P_{ik}} \\
 & { Q_i = \sum_{k \sim i} Q_{ik}}.
\end{align}
\end{subequations}
The constraints in \eqref{eqn:vol}  capture the voltage regulation goal. The  constraints in  \eqref{eqn:p} and \eqref{eqn:q}  describe the limited ability to control active/reactive power injections on each bus $i$;  $\overline{P}_i$ ($\underline{P}_i$) and $\overline{Q}_i$ ($\underline{Q}_i$),  denote the upper (lower)  limits on the amount of active and reactive power that each bus $i$ can provide, respectively. The  constraints in \eqref{eqn:pik} capture {line power flow} limits, while \eqref{eqn:lik} imposes loss limits on individual lines. 
Without loss of generality and to ease the notations in subsequent development, hereafter we assume $V_i^{ref}=1$ p.u. for all $i$. {Note that active power and reactive power need not be controllable at every bus. If for a particular bus they are not controllable, in the optimization problem we   set the bus active and/or reactive power upper and lower bounds to be equal, which essentially fixes the active and/or reactive on that bus.}

The optimization problem in \eqref{eqn:prob1} is difficult for two reasons: i) it is not convex due to the quadratic relationship between  bus voltages and powers; and  ii) depending on the size of the network, there could potentially be a large number of variables and constraints. Section \ref{sec:main} and Section \ref{sec:geometry} address i)  by convexifying the problem in  \eqref{eqn:prob1}, while Section \ref{sec:algo} addresses ii) by proposing a computationally efficient distributed algorithm to solve the resulting convexified problem.

\section{Convex Relaxation} \label{sec:main}
In this section, we state the main theoretical result, which is that under certain conditions on the angle differences between adjacent buses and the lower bounds on reactive power injections, the nonconvex problem in \eqref{eqn:prob1} can be solved exactly by solving its convex SDP relaxation. We note that SDP relaxation is not the only possible convex relaxation, e.g., \cite{Farivar12} proposes a SOCP relaxation.  In order to state the SDP relaxation, it is convenient to rewrite the problem in \eqref{eqn:prob1} in matrix form.
\subsection{Voltage Regulation Problem Formulation  in Matrix From}
Let $\bd{E}_i \in \R^{n \times n}$, with $E_{ii}=1$ and all other entries equal to zero, and define $\bd{A}_i = \frac{1}{2}(\bd{Y}^H\bd{E}_i + \bd{E}_i \bd{Y})$, and $\bd{B}_i = \frac{1}{2j}(\bd{Y}^H\bd{E}_i - \bd{E}_i \bd{Y})$.
Then, the active and reactive power injections in  bus $i$ are given by $P_i = \text{Tr}(\bd{A}_i\bd{v}\bd{v}^H)$, and $Q_i = \text{Tr}(\bd{B}_i\bd{v}\bd{v}^H)$, respectively,  where $\text{Tr}(\cdot)$ is the trace operator. For each $(i,k) \in \mathcal{E}$,  define a matrix $\bd{A}_{ik}$, with its  $(l,m)^{th}$ entry   given  by
\begin{align}
	\bd{A}_{ik}[l,m]=\left\{
	\begin{array}{ll}
		\text{Re}\{Y_{ik}\}& \text{if } l=m=i\\
		 -Y_{ik}/2 & \text{if } l=i \text{ and } m=k\\
		 -Y_{ik}^H/2 & \text{if } l=k \text{ and } m=i\\
		0 & \text{otherwise};
	\end{array} \right.
\end{align}
the active power flow through the $(i,k)$ line is given by  $P_{ik}=\Tr(\bd{A}_{ik} \bd{v} \bd{v}^H)$. Let $\bd{G}_{ik}= \bd{A}_{ik}+\bd{A}_{ki}$. Then, we can rewrite    \eqref{eqn:prob1}  as
\begin{subequations}
\label{eqn:v}
\begin{align}
\underset{\bd{v}~~}{\text{min~~}} 	& \sum_{i=1}^n{\text{Tr}(\bd{A}_i \bd{v}\bd{v}^H)}  \\
\text{s.t.~~} & |V_i| =1, \forall i   \\
 & \underline{P}_i\leq \text{Tr}(\bd{A}_i \bd{v}\bd{v}^H) \leq \overline{P}_i, \forall i 				\\
 & \underline{Q}_i\leq \text{Tr}(\bd{B}_i \bd{v}\bd{v}^H) \leq \overline{Q}_i, \forall i 		\\
 & \Tr(\bd{G}_{ik} \bd{v}\bd{v}^H) \leq \ov{L}_{ik}, \forall i \sim k \\
 & |\text{Tr}(\bd{A}_{ik} \bd{v}\bd{v}^H)| \leq \overline{P}_{ik}, \forall i\sim k.
\end{align}
\end{subequations}
Note that  the outer product $\bd{v}\bd{v}^H$ is a positive semidefinite rank-1  $n \times n$ matrix. Conversely, given a positive semidefinite rank-1 $n \times n$  matrix, it is always possible to write it as an outer product of a vector and itself. Thus, we can rewrite  \eqref{eqn:v}  as
 \vspace{-0.5mm}
\begin{subequations}
\label{eqn:Wrank}
\begin{align}
\underset{\bd{W} \sdp 0~~}{\text{min~~}} 	& \sum_{i=1}^n{\text{Tr}(\bd{A}_i \bd{W})}  \\
\text{s.t.~~} & \bd{W}[i,i] =1, \forall i  \label{eqn:vol_1} \\
 & \underline{P}_i\leq \text{Tr}(\bd{A}_i \bd{W}) \leq \overline{P}_i, \forall i 	\label{active_injection}		 \\
 & \underline{Q}_i\leq \text{Tr}(\bd{B}_i \bd{W}) \leq \overline{Q}_i, \forall i 		\\
  & \Tr(\bd{G}_{ik} \bd{W}) \leq \ov{L}_{ik}, \forall i \sim k \\
 & |\text{Tr}(\bd{A}_{ik} \bd{W})| \leq \overline{P}_{ik}, \forall i\sim k \label{thermal}  \\
 & \rank(\bd{W})=1, \label{eqn:rank}
\end{align}
\end{subequations}
where the rank-1 constraint in \eqref{eqn:rank} makes  the problem not convex.

\subsection{Convexification} \label{sec:conv}
The problem in \eqref{eqn:Wrank} is not convex due to the rank-$1$ constraint \eqref{eqn:rank};  by removing  it, we obtain a relaxation that is  convex:
\begin{subequations}
\label{eqn:W}
\begin{align}
\underset{\bd{W} \sdp 0~~}{\text{min~~}} 	& \sum_{i=1}^n{\text{Tr}(\bd{A}_i \bd{W})} \label{pdfobj2}  \\
\text{s.t.~~} & \bd{W}[i,i] =1, \forall i \label{voltagecons2}  \\
 & \underline{P}_i\leq \text{Tr}(\bd{A}_i \bd{W}) \leq \overline{P}_i, \forall i \label{realpowercons2}				 \\
 & \underline{Q}_i\leq \text{Tr}(\bd{B}_i \bd{W}) \leq \overline{Q}_i, \forall i \label{reactivepowercons2}		\\
   & \Tr(\bd{G}_{ik} \bd{W}) \leq \ov{L}_{ik}, \forall i \sim k \\
 & |\text{Tr}(\bd{A}_{ik} \bd{W})| \leq \overline{P}_{ik}, \forall i\sim k; \label{flowcons2}
\end{align}
\end{subequations}
This convex relaxation is not always tight since the rank of the solution to \eqref{eqn:W} could be greater  than $1$. Thus the solution to \eqref{eqn:W} does not always coincide with the solution to \eqref{eqn:Wrank}. However, by imposing two conditions that are widely held in practice, the non-convex problem in \eqref{eqn:Wrank} can be solved exactly by solving \eqref{eqn:W}, as stated in the following theorem.
\begin{thm}   \label{thm:main}
{Consider a power distribution  network with a tree topology.} 
Define $\te_{ik}^P$ to be the smallest positive solution to the equation $\ov{P}_{ik}=P_{ik}$, with $P_{ik}$ given in \eqref{eqn:Pik} for $|V_i|=|V_k|=1$;  and $\te_{ik}^L=\cos^{-1} (1- \frac{\ov{L}_{ik}}{2 g_{ik}})$ if $\frac{\ov{L}_{ik}}{2 g_{ik}} \leq 2$ and $\te_{ik}^L= \infty$ if $ \frac{\ov{L}_{ik}}{2 g_{ik}} >2$. Let $\ov{\te}_{ik}=\min (\te_{ik}^P,\te_{ik}^L)$. Suppose $\ov{\te}_{ik}$ satisfies
\begin{equation} \label{eqn:ang_con}
-\tan^{-1} \left( b_{ik}/g_{ik} \right) <  \ov{\te}_{ik} < \tan^{-1}\left(  b_{ik}/g_{ik} \right);
\end{equation}
and  reactive power injection lower bounds  satisfy
\begin{align} \label{second_cond}
 \underline{Q}_i <\be_i,~i=2,\dots,n,
\end{align}
 with $
 \be_i = \sum_{k: k \in \mc{C}(i)} b_{ik}-g_{ik} \sin (\tl{\te}_{ik})- b_{ik} \cos (\tl{\te}_{ik})$,
where $\mc{C}(i)$ is the set of all neighbors of $i$ and $\tilde{\te}_{ik}=\min(\tan^{-1} (\frac{g_{ik}}{b_{ik}}), \ov{\te}_{ik})$.
Let $\bd{W}^*$ be an optimal solution to the relaxed problem in \eqref{eqn:W}. Then
\begin{enumerate}
	\item If $\bd{W}^*$ is rank $1$, then $\bd{W}^*=\bd{v}^* (\bd{v}^*)^H$ for some vector $\bd{v}^*$. Furthermore, $\bd{v}^*$ is the optimal solution to the voltage regulation problem stated in \eqref{eqn:prob1}.
	\item If $\rank(\bd{W}^*) > 1$, then there is no feasible solution to the voltage regulation problem stated in \eqref{eqn:prob1}.
  \item If \eqref{eqn:W} is infeasible, then the the voltage regulation problem stated in \eqref{eqn:prob1} is infeasible.
\end{enumerate}
\end{thm}
The theorem is proved for a two-bus network in Section \ref{sec:geometry} by studying the geometry of the feasibility set of the original problem in \eqref{eqn:v} and that of its convex relaxation in \eqref{eqn:W}. The intuition and geometric insight developed by studying the two-bus network carries over to a general tree network and the full proof is provided in the Appendix.

\section{Sketch of Theorem \ref{thm:main} Proof}   \label{sec:geometry}
The insights into Theorem 1 are obtained by studying the geometry of the sets that result from the constraints on line power flows and power injections  as described in \eqref{eqn:p}--\eqref{eqn:lik}. This geometric view was explored in previous works \cite{Zhang11a,LTZ12}. Here, we revisit the results of \cite{LTZ12} and generalize them to include limits on reactive power injections.


\begin{figure}[t!]
\centering
\includegraphics[width=6.8cm]{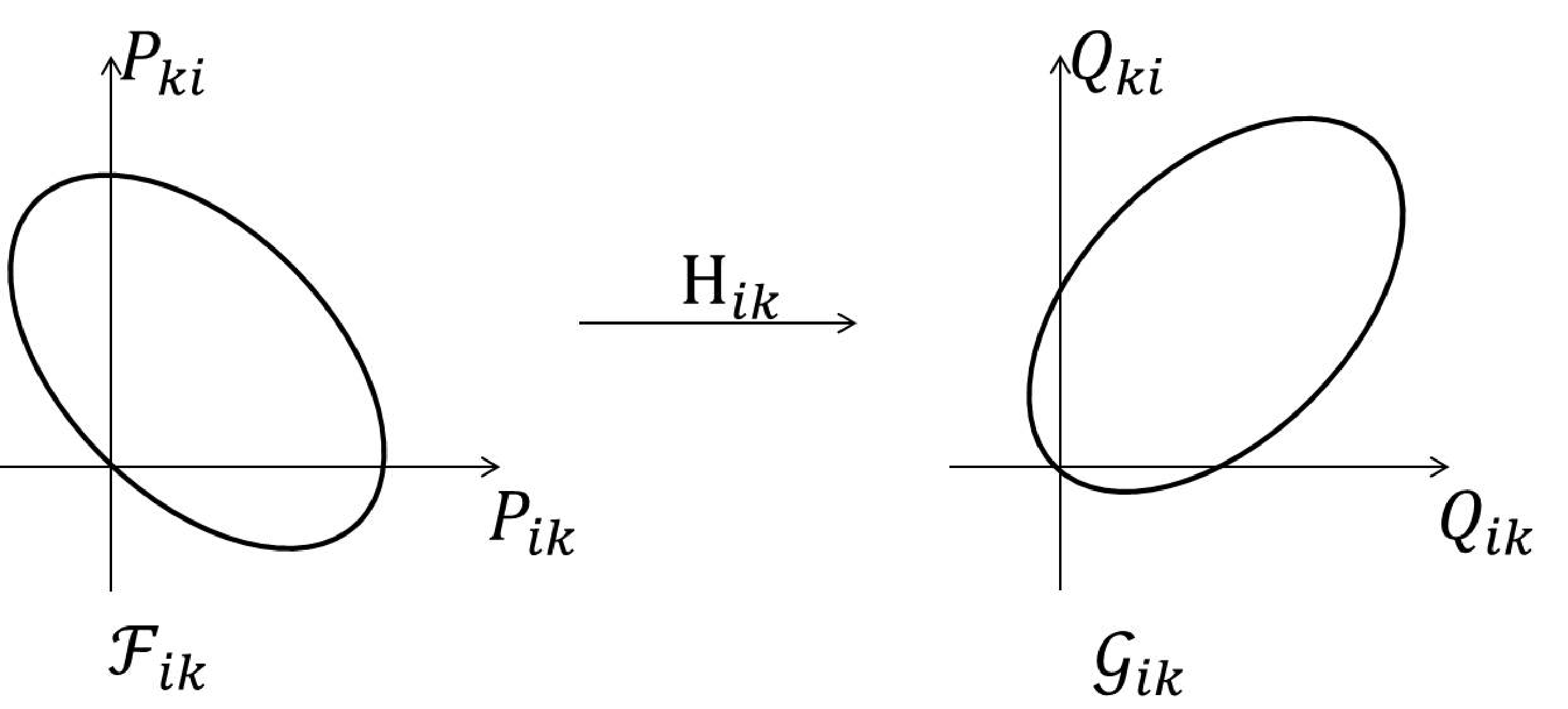}%
\caption{The active line flow region $\mc{F}_{ik}$, the reactive flow region $\mc{G}_{ik}$, and the linear transformation $\bd{H}_{ik}$ between them.}%
\label{fig:PQ}%
\end{figure}

\subsection{Active and reactive line flow regions}
First, recall from  \eqref{eqn:vol_1} that $|V_i|=|V_k|=1$ p.u. Then, let $\mc{F}_{ik} \in \mathbb{R}^2$ and $\mc{G}_{ik} \in \mathbb{R}^2$ denote the regions that contain  all the  $[P_{ik},P_{ki}]^T$ and  $[Q_{ik},Q_{ki}]^T$ that can be achieved from \eqref{eqn:Pik} and \eqref{eqn:Qik}  by varying $\te_{ik}$ between $0$ and $2\pi$; it is easy to see that for $0< \te_{ik}<2\pi$, \eqref{eqn:Pik} and \eqref{eqn:Qik} are linear transformations of a circle. Thus, as depicted in Fig.~\ref{fig:PQ},  the active and reactive line flow regions $\mc{F}_{ik}$ and $\mc{G}_{ik}$ are ellipses.  The center of  $\mc{F}_{ik}$ ($\mc{G}_{ik}$) is $[g_{ik},g_{ik}]^T$ ($[b_{ik},b_{ik}]^T$).  Its major axis is parallel to $[1,-1]^T$ ($[1,1]^T$)  and has length  $b_{ik}$ ($b_{ik}$), while its minor   axis is parallel to $[1,1]^T$ ($[1,-1]^T$)  and has length $g_{ik}$ ($g_{ik}$).  Both ellipses are related by a linear invertible mapping:  $\mc{G}_{ik}=\bd{H}_{ik} \mc{F}_{ik}$, with
\begin{align}
\bd{H}_{ik}= \frac{1}{2 b_{ik} g_{ik}}
\bma b_{ik}^2-g_{ik}^2 & b_{ik}^2 + g_{ik}^2 \\
b_{ik}^2+g_{ik}^2 & b_{ik}^2- g_{ik}  \label{transformation}
\ebma.
\end{align}

\begin{figure}[t!]
 \vspace{-0.1in}
\psfrag{Pik}{$P_{ik}$}
\psfrag{Pki}{$P_{ki}$}
\psfrag{Lik}{$\ov{L}_{ik}$}
\psfrag{oPik}{$\ov{P}_{ik}$}
\psfrag{oPki}{$\ov{P}_{ki}$}
\centering
\subfigure{
\includegraphics[width=3cm]{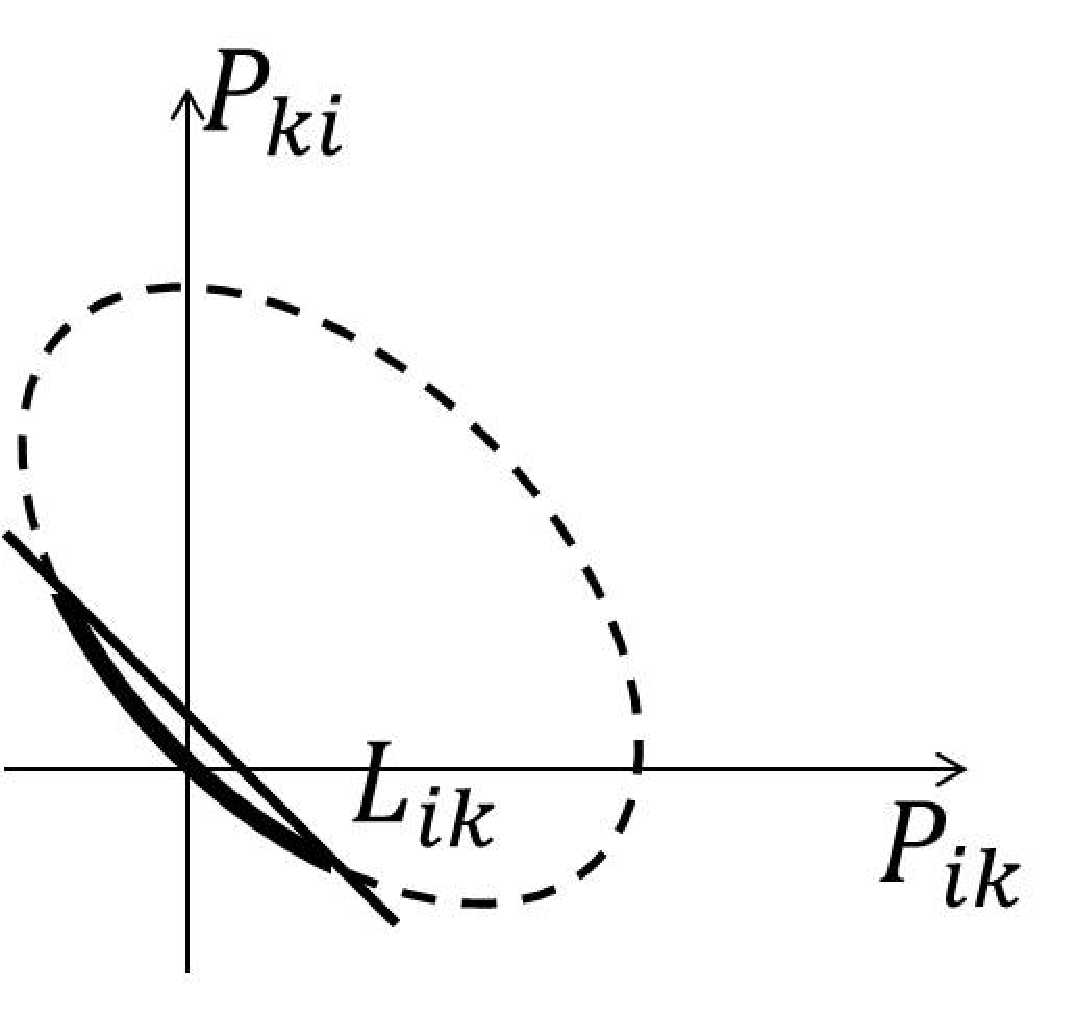}
\label{fig:Pik_loss}}
\hspace{0.5in}
\subfigure{
\includegraphics[width=3cm]{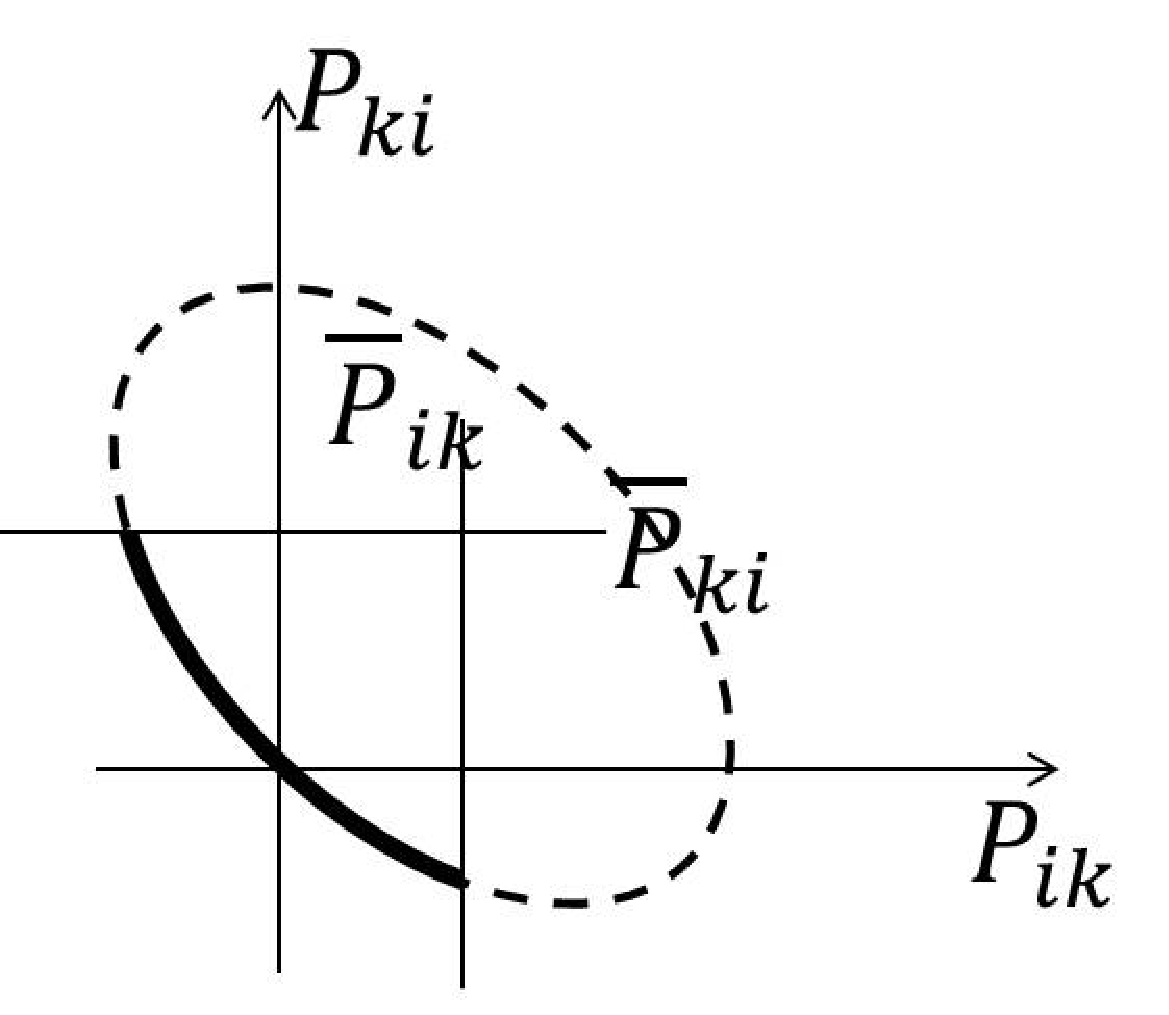}
\label{fig:Pik_line}}
\caption{The flow region under thermal loss constraints (left) and line flow constraints (right). The bold curves indicate the feasible part.}%
\label{fig:Pik}%
\vspace{-0.1in}
\end{figure}

The line flow constraints in \eqref{eqn:pik} and the thermal loss constraints in \eqref{eqn:lik}  appear  as linear constraints on the line flow regions as shown in Fig.~\ref{fig:Pik}. Thus, for each line $(i,k)$, we can replace both constraints by a single one, which has the form of the line flow constraint for properly defined upper limits. We adopt this convention in subsequent developments.
Furthermore, since  the ellipses have empty interior, this flow constraint  can be translated into angle constraints on $\te_{ik}$ of the form  $| \te_{ik}|\leq \ov{\te}_{ik}$. Conversely, an angle constraint on $\te_{ik}$ can be converted into a   flow constraint. Let $\mc{F}_{\te, ik}$ and $\mc{G}_{\te,ik}$ denote, respectively, line $(i,k)$ angle-constrained active and reactive line flow regions, then
\begin{align}
\mc{F}_{\te, ik}=\{[P_{ik},P_{ki}]^T: P_{ik}=g_{ik}[1-\cos(\theta_{ik}] + b_{ik} \sin(\theta_{ik}), \nonumber\\
P_{ki}=g_{ik}[1-\cos(\theta_{ik}] - b_{ik} \sin(\theta_{ik}),~|\te_{ik}| \leq \ov{\te}_{ik} \}, \nonumber\\
\mc{G}_{\te, ik}=\{[Q_{ik},Q_{ki}]^T: Q_{ik}=b_{ik}[1-\cos(\theta_{ik}] - g_{ik} \sin(\theta_{ik}), \nonumber\\
Q_{ki}=b_{ik}[1-\cos(\theta_{ik}] + b_{ik} \sin(\theta_{ik}),~|\te_{ik}| \leq \ov{\te}_{ik} \}. \nonumber
\end{align}

\subsection{Feasible Region of a Two-Bus Network}
\begin{figure}[t!]
  \vspace{-0.1 in}
\centering
\subfigure{
\includegraphics[width=2.9cm]{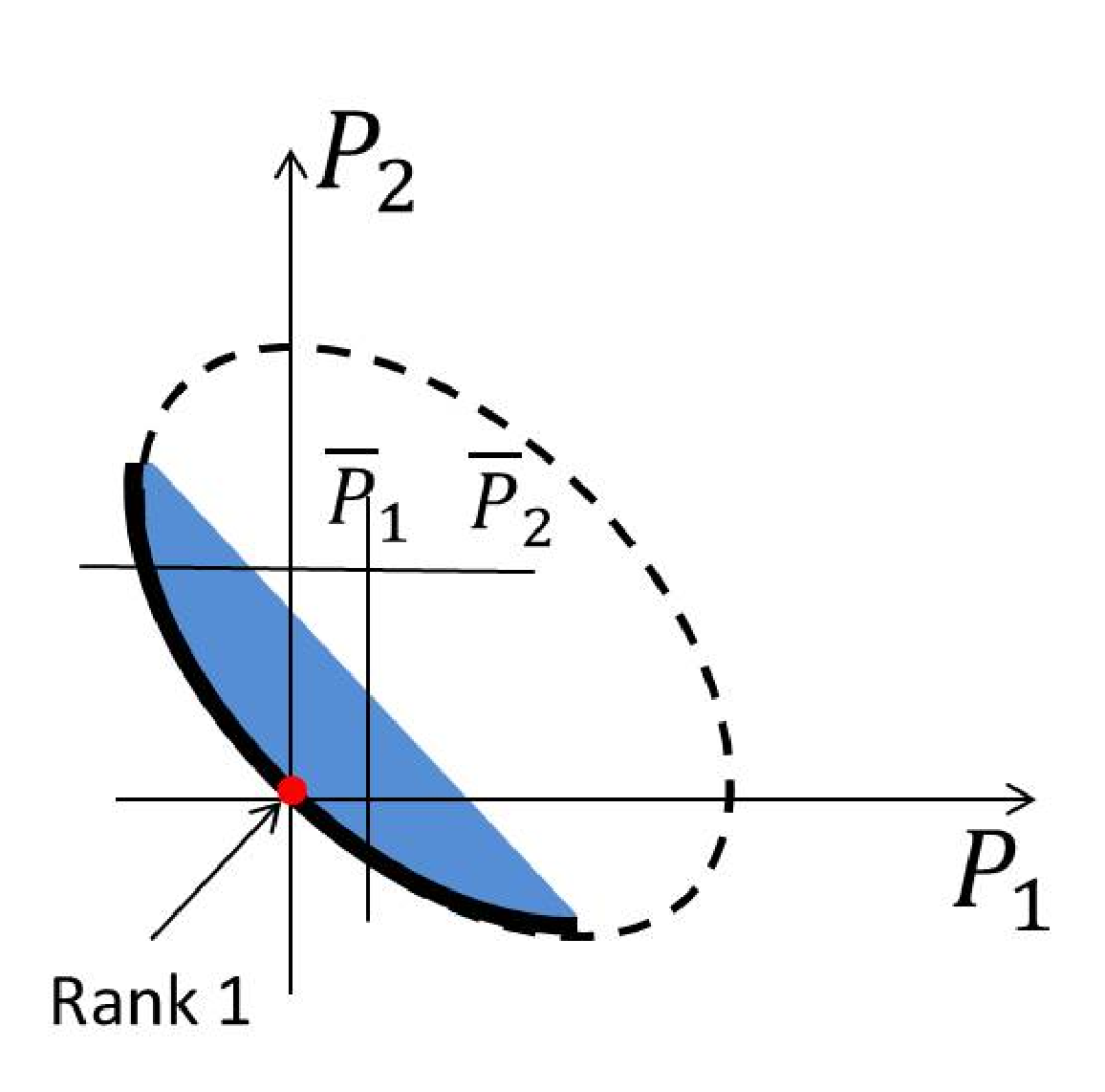}
\label{fig:2bus_P1}}
\hspace{0.75in}
\subfigure{
 \hspace{-0.1 in}
\includegraphics[width=2.6 cm]{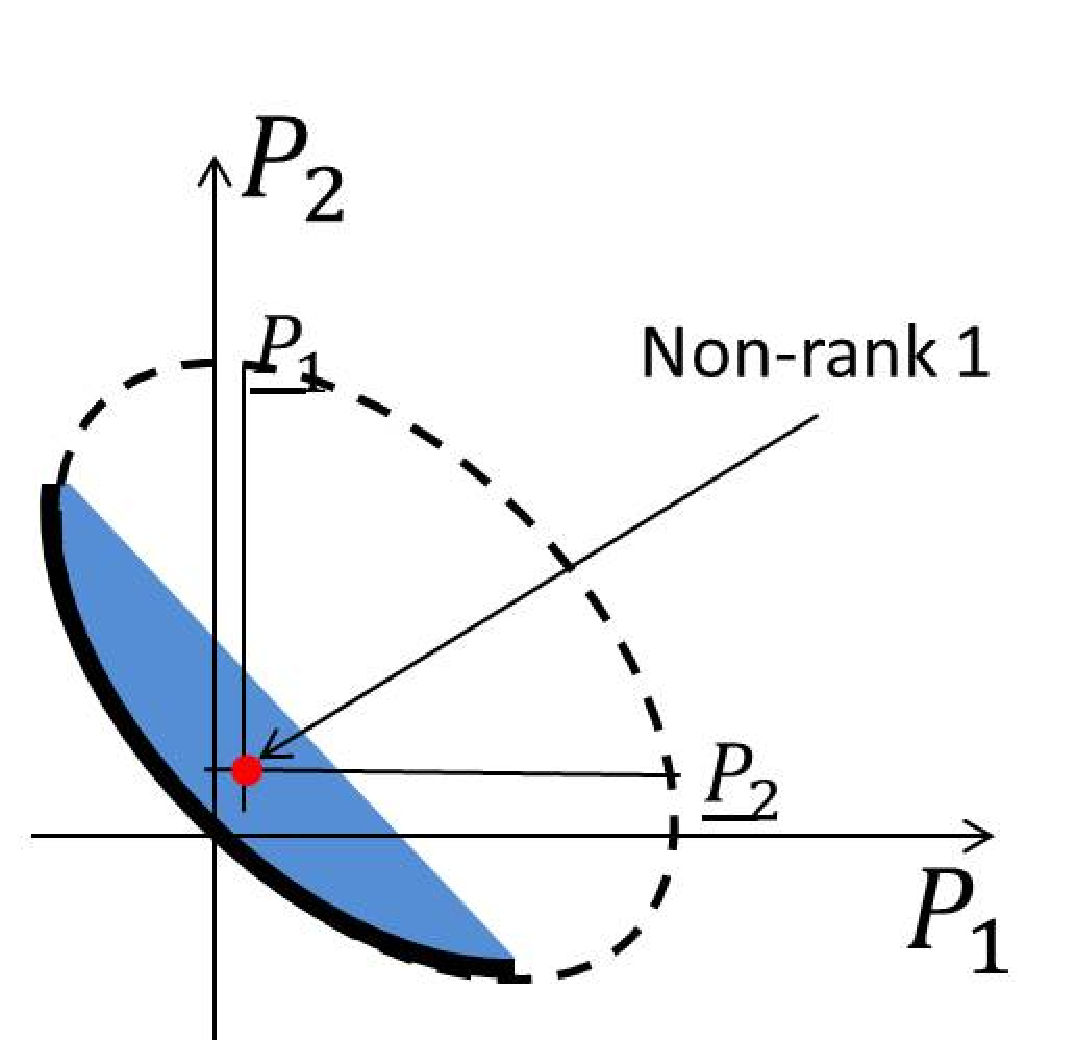}
\label{fig:2bus_P2}}
\caption{Angle-constrained active line flow region and its convex hull (filled in blue region) when relaxation is tight (left) and when relaxation does not provide solution to the original problem (right). }%
\label{fig:2bus_P}%
\vspace{-0.1in}
\end{figure}

Consider a system with only two buses connected by a line $(1,2)$, with  $P_1$ ($Q_1$) and $P_2$ ($Q_2$) denoting the active (reactive) power injections on bus $1$ and $2$, respectively; for a two-bus system, $P_1=P_{12}$ ($Q_1=Q_{12}$) and $P_2=P_{21}$ ($Q_2=Q_{21}$). The relaxed problem in \eqref{eqn:W} convexifies the feasible region of the problem \eqref{eqn:Wrank} by filling up the corresponding ellipses as shown in Fig.~\ref{fig:2bus_P}.  { Note that since the objective is to minimize the total power loss (i.e. $P_1+P_2$), the solution to the relaxed problem will  be in the lower left part of the relaxed feasible region. The relaxation is tight if the relaxed solution lies on the boundary of the ellipse, so that a rank-$1$ solution is recovered.}

The condition in \eqref{eqn:ang_con} is a constraint on the maximum angle difference across the  line. Intuitively, this angle constraint is such that only the lower left part of the line flow ellipses is feasible. For example, this condition is satisfied by the angle-constrained regions in Fig.~\ref{fig:2bus_P}; this figure shows the intersection of bus power constraints with the angle-constrained active injection regions.  In Fig.~\ref{fig:2bus_P1}, both bus power constraints are upper bounds. Since the optimal solution of the power loss problem occurs in the lower left corner, the convex relaxation is tight;  { this is an example of case $1$ in Theorem \ref{thm:main}.} In Fig.~\ref{fig:2bus_P2}, both bus power constraints are lower bounds; in this case the optimal solution   is inside the ellipses and therefore $\rank{\bd{W}^*}=2$. On the other hand, the original problem is infeasible;  { this is an example of case $2$ in Theorem \ref{thm:main}.}

It is important to note that the observations made in Fig.~\ref{fig:2bus_P} hold as long as the angle-constrained injection region only includes the lower left half of the ellipse (as described by \eqref{eqn:ang_con}). From thermal data for some common lines in \cite{Kersting06}, we expect that the angle to be constrained to $\te_{ik} \in [ -10^{\circ},10^{\circ}]$. Even for a relatively small $b_{ik}/g_{ik}$ ratio of $2$, $\ov{\te}_{ik}=\tan^{-1}(b_{ik}/g_{ik})=63.4^{\circ}$ and the condition $|\te_{ik}| < \ov{\te}_{ik}$ is always satisfied.  Therefore in most practical networks,  it is expected that the thermal constraints in the network are small enough that the condition in \eqref{eqn:ang_con} should be satisfied almost always.

\begin{figure}[t!]
\centering
\subfigure{
\includegraphics[width=3cm]{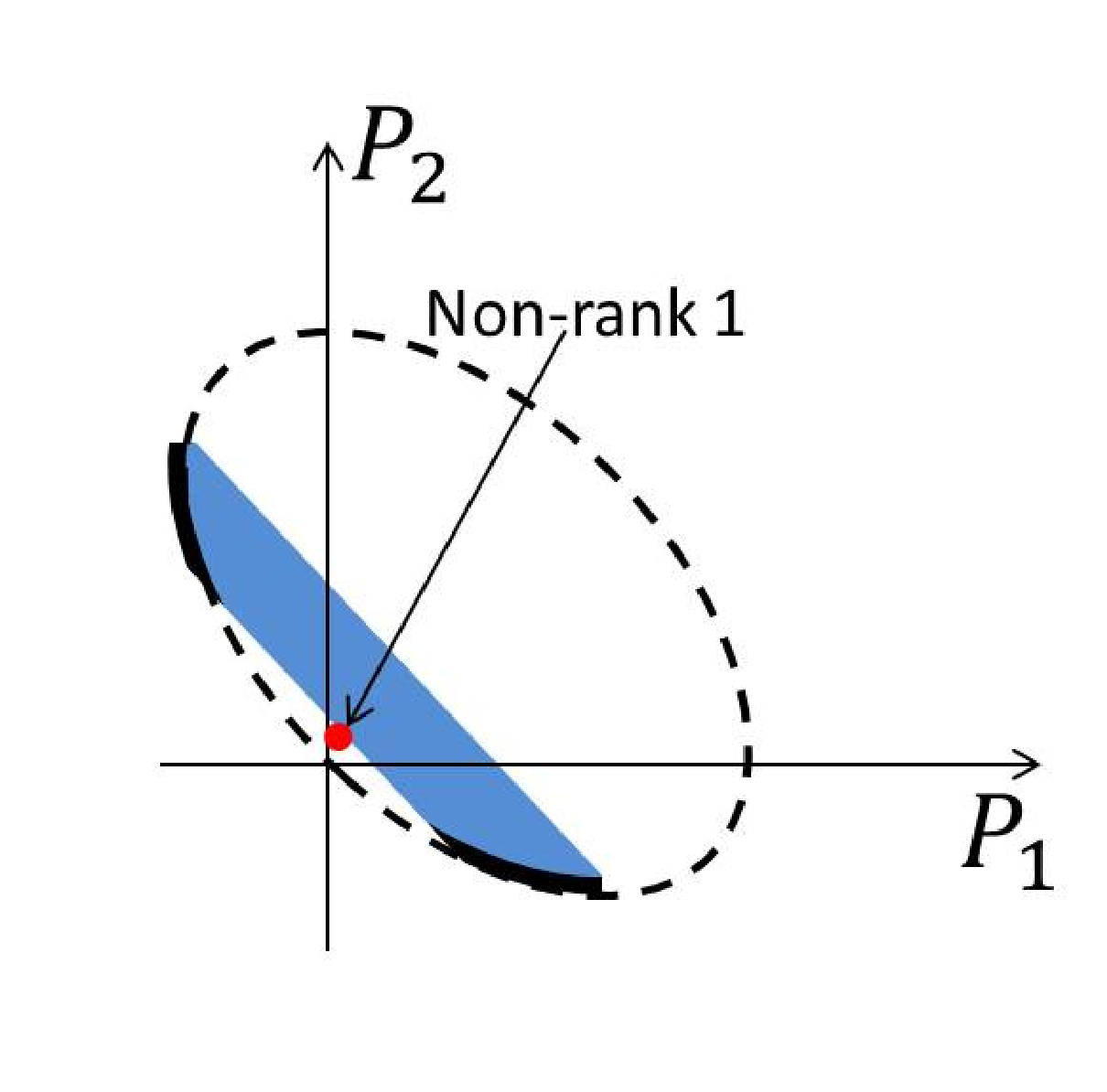}
\label{fig:PQt1}}
\hspace{0.75in}
\subfigure{
\includegraphics[width=3cm]{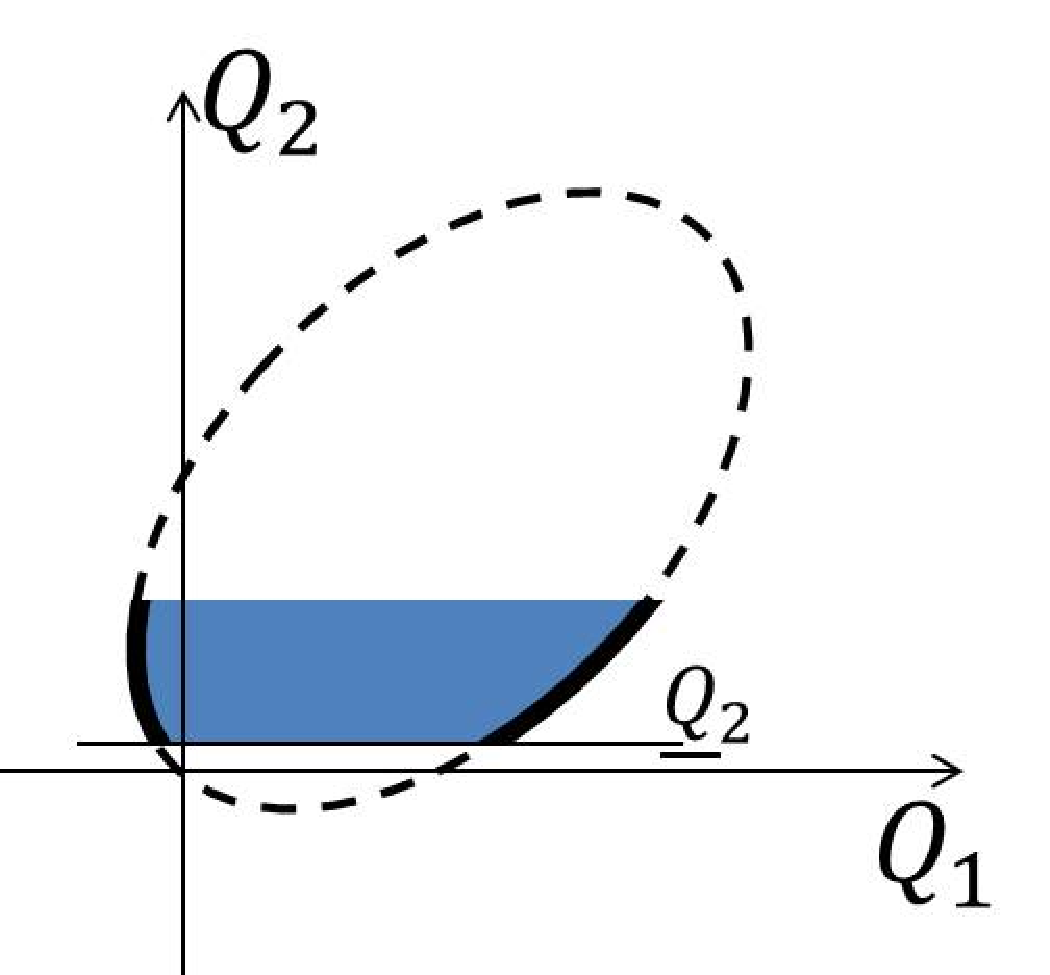}
\label{fig:PQt2}}
\caption{Active power injection region (left) and reactive power injection region (right) under reactive power injection lower bound.}%
\label{fig:PQt}%
\vspace{-0.2in}
\end{figure}

The second condition in \eqref{second_cond} is to ensure that the reactive lower bound is large enough such that $Q_i > \ul{Q}_i$ for all feasible $Q_i$. If the reactive lower bounds are tight at the optimal solution of the relaxed problem, then the rank of the optimal matrix $\bd{W}^*$ is not necessarily 1.
Figure~\ref{fig:PQt} shows the reason that the condition on the reactive power lower bounds are needed.  Figure~\ref{fig:PQt2} gives the reactive injection region with a tight reactive lower bound on bus 2. Figure~\ref{fig:PQt1} shows the corresponding active power injection region. Observe that it is possible for the optimal solution of the relaxed problem to be of rank $2$, while the original problem remains feasible. The condition $\underline{Q}_i < \be_i$ rules out this phenomenon by ensuring that the reactive power lower bounds are never tight.

\vspace{-2mm}
\subsection{General Tree Networks} \label{active_reactive_injections}
The geometrical intuition  developed for the two-bus network carries over to a general tree network due to the fact that   flows on each line are independent (no cycles), and  active and reactive power injections  can be described, respectively, as linear combinations of active and reactive line flows. These are the main ideas used in proving  Theorem \ref{thm:main}; the interested reader is referred to the Appendix for the full proof.

\vspace{-2.5mm}
\section{A Distributed Algorithm for Solving the Convexified Problem} \label{sec:algo}
In Section \ref{sec:main}, we showed that the SDP program in \eqref{eqn:W} is a convex relaxation of the voltage regulation problem in \eqref{eqn:prob1}. Since the objective is to regulate the voltages in the presence of fast-changing power injection that, e.g., arise from renewable-based generation; the optimization problem needs to be solved no slower than the time-scale at which these injections  significantly change.
General-purpose SDP solvers scale poorly as the problem size increases \cite{Boyd04}. Thus for large distribution networks with hundreds or thousands of buses, solving the SDP problem in a minute to sub-minute scale is challenging.
Furthermore standard solvers for SDP problems are centralized;  i.e., it is assumed that all the  data defining the problem is available to a single processor. However the communication infrastructure in a distribution network may not be able to transmit all the data   to a centralized location fast enough.
By exploiting the tree structure of  distribution networks, we propose a distributed algorithm to solve \eqref{eqn:W} that only requires communication between neighboring buses. This communication requirement is reasonable since that neighboring  buses are typically physically closest to each other as well. Therefore any wireless communication technology (and obviously power line communication) would enable nearest neighbors to communicate to each other. 

\subsection{Algorithm Derivation}
  The proposed algorithm consists of two stages: \textit{local optimization} and \textit{consensus}. In the local optimization stage, each node solves its own local version of the problem. In the consensus stage, neighboring nodes exchange Lagarangian multipliers obtained from the solutions to their corresponding local optimums, with the goal of equalizing the phase angle differences across a line from both of its ends.

Let $\mathcal{N}_i$ be the set of buses directly connected to bus $i$ by transmission lines, together with bus $i$ itself, i.e., $\mathcal{N}_i=\{k: k\sim i,\forall k\}\cup \{i\}$. For a $n\times n$ matrix $\bd{M}$, let $\bd{M}^{(i)}$ denote the $|\mathcal{N}_i| \times |\mathcal{N}_i|$ submatrix of $\bd{M}$ whose rows and columns are indexed according to $\mathcal{N}_i$. Similarly, for the $n\times 1$ vector $\bd{v}$, $\bd{v}^{(i)}$ is the corresponding $\mathcal{N}_i$-dimensional vector indexed by $\mathcal{N}_i$. We can rewrite  \eqref{eqn:W} as
\begin{subequations}
\label{prob3}
\begin{align}
\underset{\bd{W}^{(1)},\ldots,\bd{W}^{(n)} \sdp 0~~}{\text{min~~~~~}}  	& \sum_{i=1}^n{\text{Tr}(\bd{A}^{(i)} \bd{W}^{(i)})}   \label{pdfobj3} \\
\text{s.t.~~~~~~~~~~} & \text{diag}(\bd{W}^{(i)}) =\bd{v}^{(i)}\circ \bd{v}^{(i)}, \forall i	\label{voltagecons3}    \\
 & \underline{P}_i\leq \text{Tr}(\bd{A}^{(i)} \bd{W}^{(i)}) \leq \overline{P}_i, \forall i 			 \label{realpowercons3}	\\
 & \underline{Q}_i\leq \text{Tr}(\bd{B}^{(i)} \bd{W}^{(i)}) \leq \overline{Q}_i, \forall i 			 \label{reactivepowercons3}	\\
  &  |\text{Tr}(\bd{A}_{ik}^{(i)} \bd{W}^{(i)})| \leq \overline{P}_{ik}, \forall (i,k)\in \mathcal{E} 			 \label{flowcons3}	\\
 &W^{(i)}_{ik}=W^{(k)}_{ik}, \forall (i,k)\in \mathcal{E}, \label{lineik3} \\
 &W^{(i)}_{ki}=W^{(k)}_{ki}, \forall (i,k)\in \mathcal{E}, \label{lineki3}
\end{align}
\end{subequations}
where $\circ$ is the Hadamard product.
It is easy to verify that (\ref{pdfobj3}), (\ref{voltagecons3}), (\ref{realpowercons3}), (\ref{reactivepowercons3}), and (\ref{flowcons3}) are equivalent to  (\ref{pdfobj2}), (\ref{voltagecons2}), (\ref{realpowercons2}), (\ref{reactivepowercons2}), and (\ref{flowcons2}), respectively, as $\bd{A}_i$ in \eqref{realpowercons2}, $\bd{B}_i$ in \eqref{reactivepowercons2}, and $\bd{A}_{ik}$ in \eqref{flowcons2} have non-zero elements only at $(i,i)$, $(i,k)$, $(k,i)$, $\forall k\sim i$. Since the network is a tree, the maximal cliques are the set of adjacent nodes connected by an edge. Consequently, $\bd W^{(i)} \succeq 0$ for all $i$  is equivalent to $\bd W \succeq 0$ because the set of $\mc N_i$'s includes   all the maximal cliques of the network \cite{MCS,Lam11}. Constraints (\ref{lineik3}) and (\ref{lineki3}) are added to ensure that all $\bd{W}^{(i)}$'s  coordinate to form $\bd{W}$; in other words, $\forall (i,k)\in \mathcal{E}$, the $\theta_{ik}$'s computed from $\bd{W}^{(i)}$ and $\bd{W}^{(k)}$ should be the same.

Let $\lambda_{ik}$ be the Lagrangian multiplier of (\ref{lineik3}) for $(i,k)$ and similarly $\lambda_{ki}$ for (\ref{lineki3}). By relaxing (\ref{lineik3}) and (\ref{lineki3}), the augmented objective function is
\begin{align}
&\sum_{i=1}^n{\text{Tr}(\bd{A}^{(i)} \bd{W}^{(i)})} + \sum_{(i,k)\in\mathcal{E}}[\lambda_{ik}(W^{(i)}_{ik}-W^{(k)}_{ik}) \nn\\
&+ \lambda_{ki}(W^{(i)}_{ki}-W^{(k)}_{ki}) ] \triangleq \sum_{i=1}^n{\text{Tr}(\tilde{\bd{A}}^{(i)} \bd{W}^{(i)})}, \label{augmentedobj}
\end{align}
where $\tilde{\bd{A}}^{(i)}$ is also Hermitian, and   its $(i,k)^{th}$  entry is i) $\tilde{A}^{(i)}_{ik}=A^{(i)}_{ik}$ if $i=k$, $\tilde{A}^{(i)}_{ik}=A^{(i)}_{ik}+\lambda_{ik}^H$ if $i<k$, and iii) $\tilde{A}^{(i)}_{ik}=A^{(i)}_{ik}-\lambda_{ik}^H$ if $i>k$.
With (\ref{augmentedobj}), problem (\ref{prob3}) can be divided into $n$ separable subproblems and the $i$th subproblem corresponds to bus $i$, defined as follows:
\begin{subequations}
\label{prob4}
\begin{align}
\underset{\bd{W}^{(i)} \sdp 0~~}{\text{min~~~}}	& \text{Tr}(\tilde{\bd{A}}^{(i)}\bd{W}^{(i)})   \label{pdfobj4} \\
\text{s.t.~~~~} & \text{diag}(\bd{W}^{(i)}) =\bd{v}^{(i)}\circ \bd{v}^{(i)}	\label{voltagecons4}    \\
 & \underline{P}_i\leq \text{Tr}(\bd{A}^{(i)} \bd{W}^{(i)}) \leq \overline{P}_i 			\label{realpowercons4}	 \\
 & \underline{Q}_i\leq \text{Tr}(\bd{B}^{(i)} \bd{W}^{(i)}) \leq \overline{Q}_i 			\label{reactivepowercons4}	 \\
  &  |\text{Tr}(\bd{A}_{ik}^{(i)} \bd{W}^{(i)})| \leq \overline{P}_{ik} ,\quad\forall k\sim i.			 \label{flowcons4}	
\end{align}
\end{subequations}
We denote the feasible region described by \eqref{voltagecons4}--\eqref{flowcons4} together with $\bd{W}^{(i)}\sdp 0$ of Subproblem $i$ by $\mathcal{C}_i$. Define $g_i(\lambda_{ik})\triangleq\inf_{\bd{W}^{(i)}\in \mathcal{C}_i}\{\text{Tr}(\tilde{\bd{A}}^{(i)}\bd{W}^{(i)})\}.$
 The gradient of $-g_i$
  at $\lambda_{ik}$ is $W^{(i)*}_{ik}$, which is the $(i,k)$th element of the optimal $\bd{W}^{(i)*}$ of $g_i$ determined by solving the $i$th subproblem \eqref{prob4}. Similarly, that of $-g_k$ at $\lambda_{ik}$ is $-W^{(k)*}_{ik}$. Therefore, the gradient of $-(g_i+g_k)$ is then $W^{(i)*}_{ik}-W^{(k)*}_{ik}$. Let $W^{(i)}_{ik}[t]$ and $W^{(k)}_{ik}[t]$ be $W^{(i)*}_{ik}$ and $W^{(k)*}_{ik}$ determined at time $t$, respectively. By gradient ascent, at time $t+1$, we update $\lambda_{ik}$ by
\begin{align}
\lambda_{ik}[t+1] = \lambda_{ik}[t] + \alpha[t](W^{(i)}_{ik}[t]-W^{(k)}_{ik}[t]), \label{lambdaupdate}
\end{align}
where $\alpha[t]>0$ and $\lambda_{ik}[t]$ are the step size and $\lambda_{ik}$ at time $t$, respectively. The value of $\lambda_{ki}[t+1]$ can be directly computed from $\lambda_{ik}[t+1]$ as $\lambda_{ki} = \lambda_{ki}^H$.
The Lagrangian multiplier $\lambda_{ik}$ is only defined for the line $(i,k)$ and the two buses at the ends of the edge, i.e., buses $i$ and $k$, are required to manipulate  $\lambda_{ik}$. The purpose of  (\ref{lambdaupdate}) is to make $W^{(i)}_{ik}$ and $W^{(k)}_{ik}$ as close to each other as possible with the help of $\lambda_{ik}$. The iteration in  (\ref{lambdaupdate}) can be computed either by bus $i$ or by bus $k$ and it is independent of all other buses and edges. Whenever both the $i$th and $k$th subproblems have been computed and so $W^{(i)}_{ik}$ and $W^{(k)}_{ik}$ have been updated, then $\lambda_{ik}$ can then be updated by using (\ref{lambdaupdate}).


 The optimization problem comprised of \eqref{augmentedobj}, together with all the constraints \eqref{voltagecons4}--\eqref{flowcons4}, imposed on the subproblems, is a dual problem of \eqref{prob3}.
When all $\lambda_{ik}$'s are optimal, $W^{(i)}_{ik}$ will be equal to $W^{(k)}_{ik}$ for all $(i,k)$'s and thus the duality gap is zero. Accordingly, we can construct the optimal $\bd{W}^*$ of problem (\ref{eqn:rank}) from the values of the $W^{(k)}_{ik}$'s.  Algorithm \ref{alg:DA} can be seen as a dual decomposition algorithm, where the constraints on the consistence of line flows are dualized. Due to the convexity of \eqref{prob3}, Algorithm \ref{alg:DA} converges to the optimal solution \cite{Terelius11}. 
{The iterative algorithm \eqref{lambdaupdate} is a subgradient method and several step size rules can be applied to specify $\alpha[t]$, e.g., constant step size,  and non-summable diminishing step size  $\alpha[t]=a/\sqrt{t}$, where $a>0$ \cite{Boyd08}.}


\begin{algorithm}[!h]
\caption{{ Distributed Algorithm} \label{algorithm_1}}
\label{alg:DA}
\footnotesize
\begin{tabbing}
Given a $n$-bus network\\
1. \=\textbf{while} $|W^{(i)}_{ik}-W^{(k)}_{ik}|>\delta$ for any $(i,k)\in \mathcal{E}$ \textbf{do}\\
\>2. \textbf{for} \= each bus $i$ (in parallel) \textbf{do}\\
\>\>3. Given $\lambda_{ik},\forall k\sim i$, solve (\ref{prob4})\\
\>\>4. Return $W^{(i)}_{ik}, \forall k$\\
\>5. \textbf{end for}\\
\>6.\= Given $W^{(i)}_{ik}$ and $W^{(k)}_{ik}$, update $\lambda_{ik}$ with (\ref{lambdaupdate}) (in parallel)\\
7. \textbf{end while}
\end{tabbing}
\vspace{-1mm}
\end{algorithm}

\vspace{-3.5mm}
\subsection{Feasibility}
When the buses determine their own limits on active and reactive powers independently, an infeasible problem might result, i.e., an empty feasible region. When there exists a central authority having all the bus power information, we can check the feasibility easily. Otherwise, it is necessary for the buses to declare infeasibility.

One sufficient condition for infeasibility of the the problem is that there exists an infeasible subproblem (\ref{prob4}) for any bus.
If any bus finds an infeasible subproblem, it is sufficient to say that the whole problem is infeasible. To proceed further, the bus with an infeasible subproblem should adjust its own active and reactive power limits so as to make the subproblem feasible.
A necessary and sufficient condition for infeasibility is that  $W^{(i)}_{ik}$ and $W^{(k)}_{ik}$ never match for some $(i,k)\in \mathcal{E}$ when Algorithm \ref{alg:DA} evolves. If this happens on edge $(i,k)$, either bus $i$ or bus $k$ or both constitute the infeasibility.

\begin{figure}[!t]
  \centering
  \subfigure[Network structure.]{
    \includegraphics[width=4.2cm]{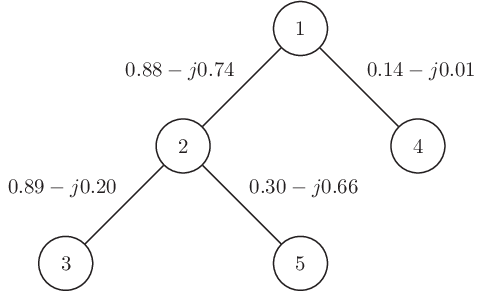}
    \label{fig:toy1}}
  \subfigure[Convergence curve.]{
  \includegraphics[width=4cm]{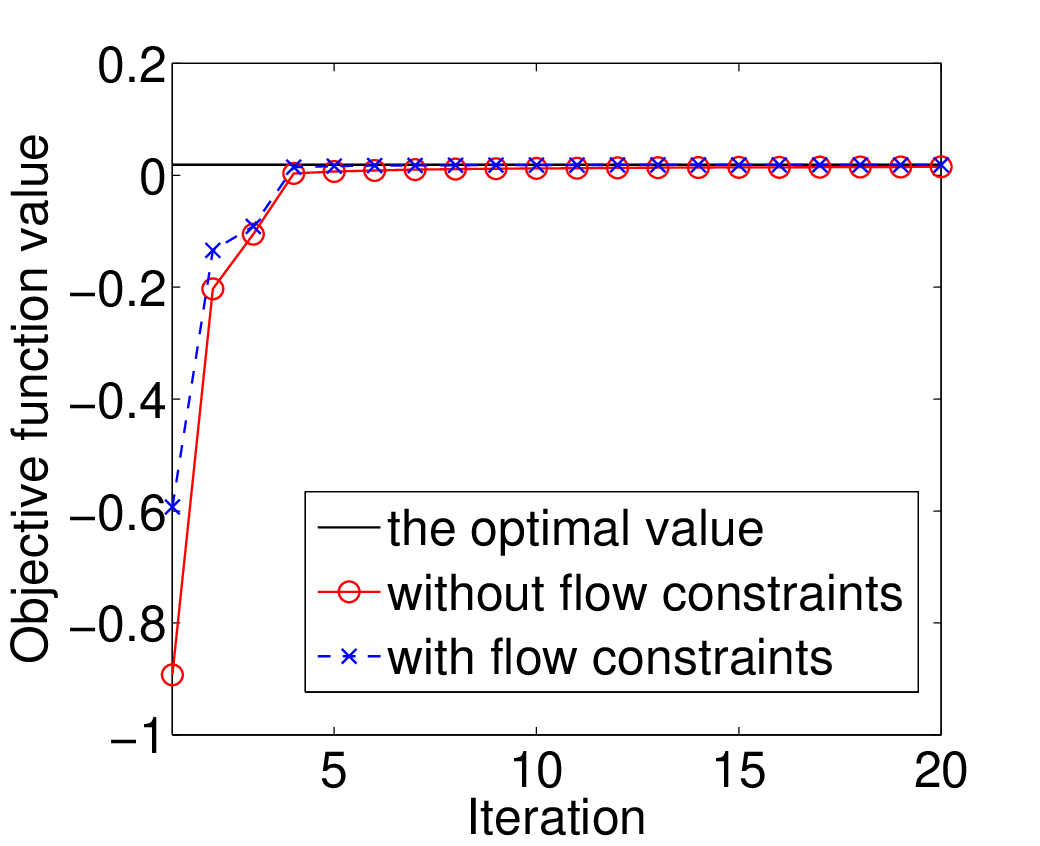}
    \label{fig:toy1lb}}
    \caption{A 5-bus example.}
\end{figure}

\vspace{-3mm}
\subsection{Numerical Performance Enhancements} \label{enhancements}
Consider the five-bus network given in Fig.~\ref{fig:toy1}. Assuming that all $\lambda_{ik}$'s are updated at the end of each iteration, the progress of Algorithm \ref{alg:DA} (the curve without power flow constraints) and the target optimal objective value are shown in Fig.~\ref{fig:toy1lb}. At iteration 20, when we sum the objective function values of all the subproblems, the sum still has around $20\%$ difference to the optimal one. Even for a small network, it may take a long time for the algorithm to converge to the global optimal solution.
Next, we provide some enhancements that improve the algorithm convergence speed.

\begin{table}[!t]
\label{table:toy1}
\caption{Bus information of the five-bus example}
\scriptsize
\centering
\begin{tabular}{c c c c c c}
  \hline\hline
Bus 	& $\overline{P}$ & $\underline{P}$ & $\overline{Q}$ & $\underline{Q}$ & $\overline{V}$\\
\hline
1	& 5.2844	& -5.4692	& 5.5798	& -5.7604	& 1.2247\\
2	& -0.0648	& -0.0988 	& 0.5298	& 0		& 1.1509\\
3	& -0.0423	& -0.5828	& 0.6405	& 0		& 1.1103\\
4	& -0.0334	& -0.5155	& 0.2091	& 0		& 0.9762\\
5	& -0.0226 & -0.4329	& 0.3798	& 0		& 1.1400\\
\hline\hline
\end{tabular}
\vspace{-0.1in}
\end{table}

\subsubsection{Power Flow Constraints}
Constraint \eqref{flowcons4} means that the active power can flow in any direction on the edge $(i,k)$ as long as its magnitude does not exceed the limit $\overline{P}_{ik}$. Assume that the global optimal solution $\bd{W}^*$ exists.
Our decomposition allows us to compute $W_{ik}^*$ separately by buses $i$ and $k$, in which each bus determines its local version of $W_{ik}^*$, e.g., $W^{(i)}_{ik}$ for bus $i$. Then \eqref{lambdaupdate} brings both $W^{(i)}_{ik}$ and $W^{(k)}_{ik}$ towards $W_{ik}^*$ by just equalizing $W^{(i)}_{ik}$ and $W^{(k)}_{ik}$. If the feasible regions $\mathcal{C}_i$ and $\mathcal{C}_k$
are smaller, it will be easier for \eqref{lambdaupdate} to reduce the discrepancy between $W^{(i)}_{ik}$ and $W^{(k)}_{ik}$.

The additional assumption we make is that all buses are net consumers of active power except the feeder; that is, $P_i \leq 0$ for $i=2,3,\dots,n$.
For faster convergence rate, we assume that active power flows from buses $i$ to $k$ along the edge $(i,k)$ with $i<k$, i.e., $P_{ik}\geq 0$. Note this assumption is not necessary for the theoretical results in Section III, but it makes the algorithm much simpler. In practice, DERs are currently not allowed to cause reverse current flow due to protection issues, but it would be interesting to generalize our  algorithm to also handle this case.  With this assumption, we can re-write \eqref{flowcons4} as
\begin{align}
0\leq P_{ik} = \text{Tr}(\bd{A}_{ik}^{(i)}\bd{W}^{(i)})\leq \overline{P}_{ik}, \label{improved_flow_limit_i}\\
-\overline{P}_{ik} \leq P_{ki} = \text{Tr}(\bd{A}_{ik}^{(k)}\bd{W}^{(k)})\leq 0, \label{improved_flow_limit_k}
\end{align}
from the perspectives of buses $i$ and $k$, respectively.
We can actually replace \eqref{flowcons4} for $(i,k)$ of Subproblem $i$ by \eqref{improved_flow_limit_i} and similarly \eqref{flowcons4} for $(i,k)$ of Subproblem $k$ by \eqref{improved_flow_limit_k}. If we apply the same logic to all edges connecting to bus $i$, we can construct a smaller feasible region $\hat{\mathcal{C}}_i$ for Subproblem $i$. For the edge $(i,k)$, the constructions of $\hat{\mathcal{C}}_i$ and $\hat{\mathcal{C}}_k$ can help $W^{(i)}_{ik}$ and $W^{(k)}_{ik}$ converge to $W_{ik}^*$ faster.

 \vspace{2mm}
With this modification, the progress of the algorithm for  the five-bus example is also depicted in Fig.~\ref{fig:toy1lb}, where  we can see that the algorithm converges faster.

\subsubsection{Feasible Solution Generation}

When the algorithm converges,
we  have that
\begin{align}
\text{Tr}(\bd{A}^{(i)} \bd{W}^{(i)})=\text{Tr}(\tilde{\bd{A}}^{(i)} \bd{W}^{(i)}), \quad\forall i, \label{optimalobj}
\end{align}
which holds when all its associated $\lambda_{ik}$'s are optimal; this is equivalent to have both of the following held:
\begin{align}
\text{Tr}(\bd{A}^{(i)} \bd{W}^{(i)}) = \text{Tr}(\bd{A}^{(i)} \bd{W}^{(i)*}) \Leftrightarrow P_i = P_i^*, \quad \forall i, \label{optimalP}\\
\text{Tr}(\bd{B}^{(i)} \bd{W}^{(i)}) = \text{Tr}(\bd{B}^{(i)} \bd{W}^{(i)*}) \Leftrightarrow Q_i = Q_i^*, \quad \forall i. \label{optimalQ}
\end{align}
In other words, Algorithm \ref{alg:DA} tries to find the the optimal active and reactive power pair $[P_i^*,Q_i^*]^T$ for each bus $i$ by manipulating $\lambda_{ik}$'s defined for the corresponding lines. The more lines are connected to a bus (i.e., the more $\lambda_{ik}$'s  it involves), the more difficult is for  (\ref{optimalP}) and (\ref{optimalQ}) to hold. The $[P_i,Q_i]^T$ pair affects the $[P_k,Q_k]^T$ pair through $\lambda_{ik}$. Consider the situation where edge $(i,k)$ is the only line connected to bus $k$ except for  bus $i$. When $[P_k,Q_k]^T$ becomes optimal, this helps bus $i$ converge in the sense that this reduces the variations of $[P_i,Q_i]^T$ induced from bus $k$. When Algorithm \ref{alg:DA} evolves, the $[P_k,Q_k]^T$ of leaf bus $k$ converges first as a leaf bus has only one edge. Then, we have the buses connected to the leaf buses converged. We continue this process and finally go up to  the feeder.

For any leaf node $k$, we have $P_k = P_{ki}$ and $Q_k = Q_{ki}$, where bus $i$ is the only bus connected to bus $k$. When the algorithm evolves,
we obtain $W^{(k)*}_{ik}$ from the solution of the $k$th subproblem (\ref{prob4})
when $\text{Tr}(\bd{A}^{(k)} \bd{W}^{(k)})$ and $\text{Tr}(\bd{B}^{(k)} \bd{W}^{(k)})$ are equal to $P_k^*$ and $Q_k^*$, respectively. Once we have fixed  $W^{(k)*}_{ik}$, we can add the constraint $W^{(i)}_{ik}=W^{(k)*}_{ik}$ to the $i$th subproblem for bus $i$ by passing a message containing the value of $W^{(k)*}_{ik}$ from bus $k$ to bus $i$.
In matrix form, this constraint is equivalent to $\text{Tr}(\bd{C}^{(i)}\bd{W}^{(i)})=\text{Re}\{W^{(k)*}_{ik}\}$ and $\text{Tr}(\bd{D}^{(i)}\bd{W}^{(i)})=\text{Im}\{W^{(k)*}_{ik}\}$,
where $\bd{C}^{(i)}=(C^{(i)}_{lm}, l,m\in \mathcal{N}_i)$, with $C^{(i)}_{lm}=\frac{1}{2}$ if $ l=i \text{ and } m=k$, $C^{(i)}_{lm}=\frac{1}{2}$ if $l=k \text{ and } m=i$, and $C^{(i)}_{lm}=0$ otherwise;
and
$\bd{D}^{(i)}=(D^{(i)}_{lm}, l,m\in \mathcal{N}_i)$, with $D^{(i)}_{lm}=\frac{1}{2}j$ if $l=i \text{ and } m=k$, $D^{(i)}_{lm}=-\frac{1}{2}j $ if $ l=k \text{ and } m=i$, and $D^{(i)}_{lm}=0$ otherwise.
In this case, we reduce the $n$-bus network into the $(n-1)$-bus one by removing bus $k$. When all other buses with positive active power flown from bus $i$ (i.e. $\{l: l\sim i, l>i\}$) have been fixed and ``removed'', bus $i$ becomes a leaf bus in the reduced network. This process continues until we find all $W_{lm}^*, \forall (l,m)\in \mathcal{E}$. The global solution $\bd{W}^*$ can be constructed from those $W_{lm}^*$'s. However, for any bus $k$, if we fix $P_k$ and $Q_k$ which is not optimal, these errors will make its connecting bus $i$ being fixed afterwards result in incorrect $P_i$ and $Q_i$, which are not optimal either.  To achieve this, we
 observe $P_i[t]$ and $Q_i[t]$ for a certain time period and check if their variations are significant. Assume that we are at time $t$, for the active power, we can keep track of the previous $T$ $P_i$'s and the current $P_i[t]$, i.e. $[P_i [t-T],P_i[t-T+1],\ldots,P_i[t]]^T$. We can say that $P_i[t]$ has converged if its cumulative change is less than a certain threshold $\gamma$ (e.g. $10^{-4}$), i.e.,
\begin{align}
\sum_{k=0}^{T-1}\frac{|P_i[t-T+k] - P_i[t-T+k+1]|}{|P_i[t-T+k]|}<\gamma; \label{P_convergence}
\end{align}
with a similar condition for the reactive power.

\begin{figure}[t!]
\centering
\vspace{-0.1in}
\includegraphics[scale=0.7]{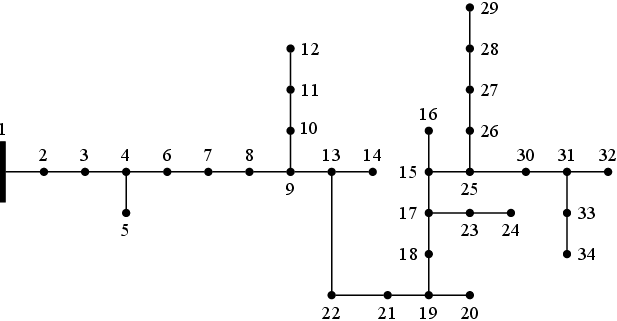}
\caption{34-bus system: electrical network graph. There are tap changing transformers between buses 7 and 8, buses 17 and 18, and buses 17 and 23.}
\label{fig:feeder34}
\vspace{-0.4in}
\end{figure}

\subsubsection{Hot Start}
The problem needs to be solved repeatedly; when there are changes to the active/reactive limits at any bus, we apply Algorithm \ref{alg:DA} to the problem again.
In each update, we usually have small variation between the new $\overline{P}_i$ and the previous ones and also for $\underline{P}_i$.
Thus, in subsequent instances of the problem, the optimal angle difference across each line usually does not vary significantly.
Therefore, we can set $\lambda_{ik}[0]$ with the optimal $\lambda_{ik}^*$ which can be determined from the previous optimal $W_{ik}^*$.

\section{Case Studies} \label{sec:simulation}

\begin{figure}[t!]
 \vspace{-0.1in}
\centering
\includegraphics[scale=0.25]{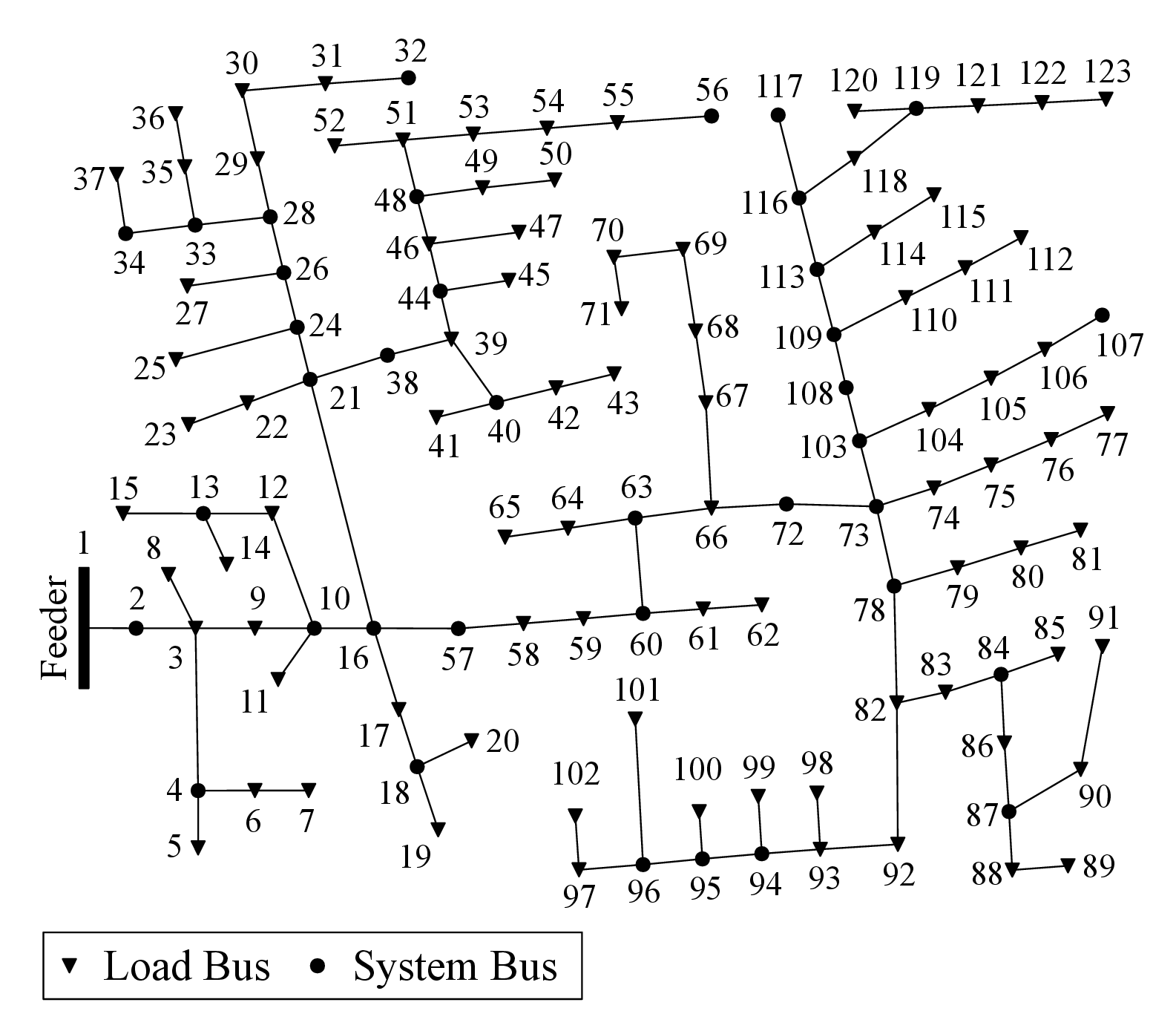}
\caption{{ 123-bus system: electrical network graph. There are tap changing transformers between buses 12 and 13, buses 28 and 33, and buses 72 and 73. Every load bus is assumed to have some capability to provide reactive power, proportional to their active power demands.}}
 \vspace{-0.2in}
\label{fig:123bus}
\end{figure}

We  test  the performance of Algorithm~\ref{algorithm_1} on the  IEEE 34- and 123-bus test systems  \cite{feeder}; the data for these systems can be found in  \cite{testfeeders}. The topology for the 34-bus system is displayed in Fig.~\ref{fig:feeder34}, while the topology for the 123-bus system is displayed in Fig.~\ref{fig:123bus}. {All simulations were performed on a MacBookPro6,2, and each one was terminated when 300 iterations were reached.}

\begin{figure}[b!]
     \vspace{-0.2in}
  \centering
  \subfigure[Entire Time Horizon. \label{entire_horizon}]{
    \includegraphics[height=3.4cm]{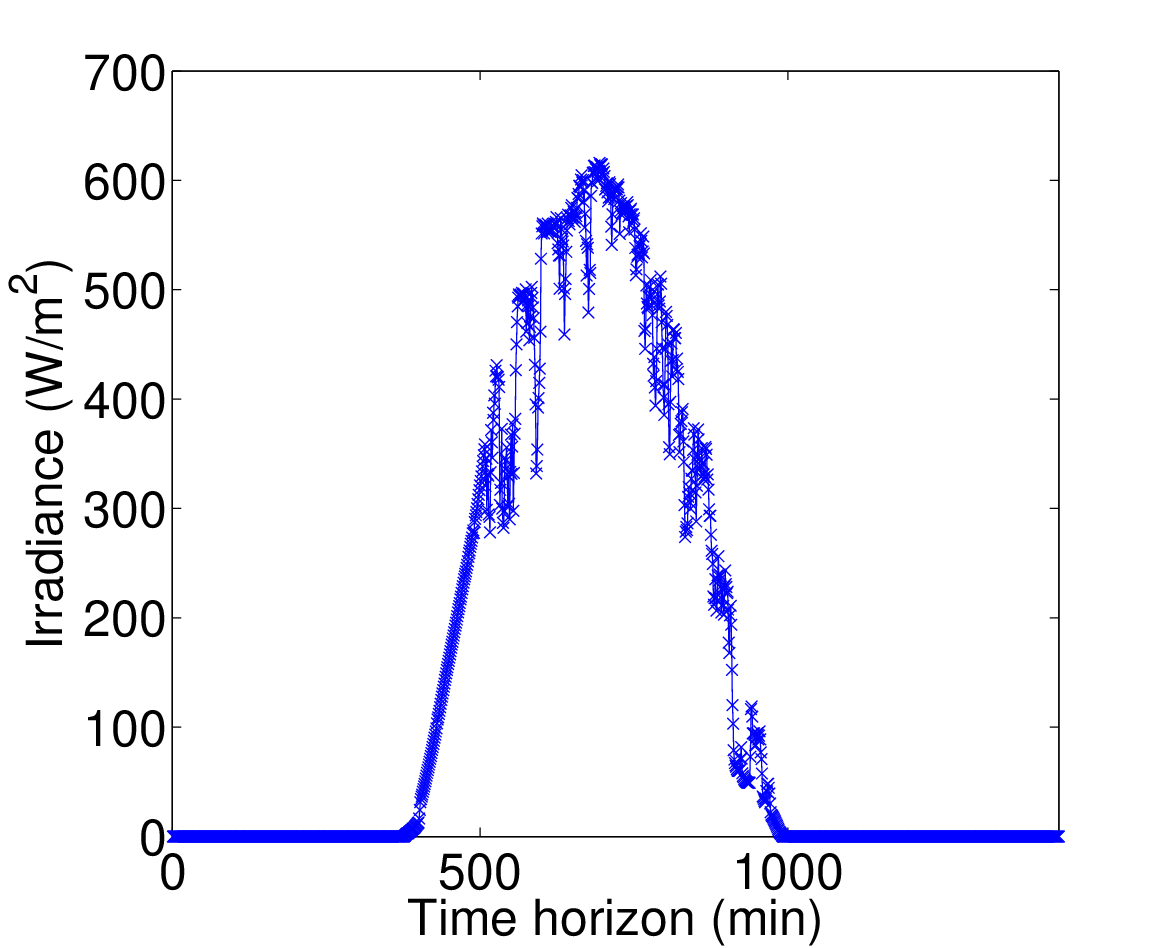}}
    \hspace{-0.2in}
   \subfigure[Minute 781 to 840. \label{hour_horizon}]{
   \includegraphics[height=3.4cm]{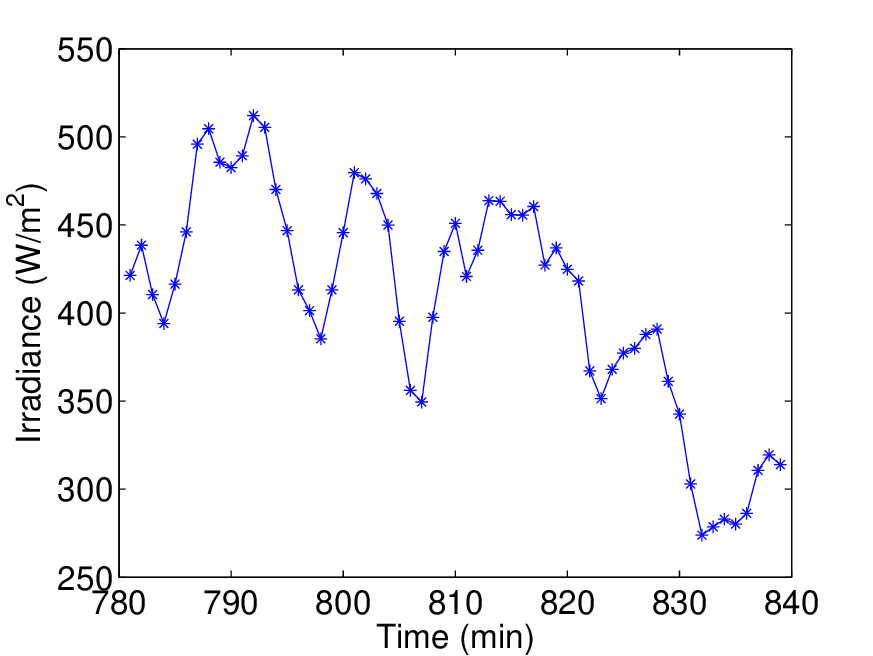}
   \label{fig:pvdata_2}}
    \caption{Irradiance of a particular day in November 2011 \cite{nrel}. \label{fig:pvdata}}
      \vspace{-0.1in}
\end{figure}

Assume that, for both test systems, the nominal load on each bus $i$, denoted  by $\hat{P}_i$, is  specified by the datasets in \cite{testfeeders}.  Additionally, we assume that connected to  each bus $i$, there are  energy storage devices  and  PV-based electricity generation resources,  which can supply active power, denoted by $P^{PV}_i$, to the bus locally, i.e., their net effect is to reduce the load. If all $P^{PV}_i$ is consumed locally, then the active power injection at bus $i$ will be $\overline{P}_i =\hat{P}_i+P^{PV}_i\leq 0$. The computed optimal $P_i^*\in [\hat{P}_i,\overline{P}_i],~i=2,\ldots,n,$ will then be adjusted by controlling the amount of power from the PV devices which will be stored at the local storage device.
Let $\hat{Q}_i$ be the nominal reactive power injection at bus $i$. By following  \cite{Zou12}, the power electronics interface  of the PV installations is assumed to be able to supply reactive power in a range that is sufficient to cancel the nominal reactive power \cite{Zou12}. Therefore, we assume that the reactive power can be adjusted in the ranges specified by i) $Q_i \in [0,1.2\hat{Q}_i]$, if $\hat{Q}_i\geq 0$, and ii) $Q_i=[-1.2\hat{Q}_i,0]$ otherwise.

We consider the one-minute resolution irradiance data in Fig.~\ref{entire_horizon}, which correspond to  a particular day in November 2011 collected at the University of Nevada \cite{nrel}; the $P^{PV}_i$'s vary in accordance  to the variation of this irradiance data. Assume that the PV systems connected to  bus $i$ can provide up to $20\%$ of  the nominal load $\hat{P}_i$ at that bus. Thus, the maximum $P^{PV}_i$, which is proportional to the respective  $\hat{P}_i$, is different for different buses. As it can be seen in  Fig.~\ref{fig:pvdata}, since there is only radiation between the $377^{th}$ and $991^{th}$ minutes, for all numerical examples, we define a time horizon of $[377,991]$, and execute  Algorithm~\ref{algorithm_1}  every minute within this time horizon. Recall that   Algorithm~\ref{algorithm_1}  requires inputs of Lagrangian multipliers as the starting points. In minute $t$, where $t\in[377,991]$, the inputs to Algorithm~\ref{algorithm_1} are the Lagrangian multipliers computed by Algorithm~\ref{algorithm_1} at time $t-1$. Moreover each Lagrangian multiplier is only stored and manipulated by the two buses at the two ends of the corresponding transmission line. Initially, i.e., at $t=377$, the  Lagrangian multipliers are computed from the nominal system settings. 
{At each step, we check if Algorithm 1 converges with a stopping criterion. It is deemed converged if the Euclidean norm of the change of the Lagrangian multiplier is smaller than some tolerance, i.e., $\|\lambda[t+1]-\lambda[t]\|<1e-6$.} 

\begin{figure}[t!]
\centering
\subfigure[Active power injections. ]{
\includegraphics[scale=0.25]{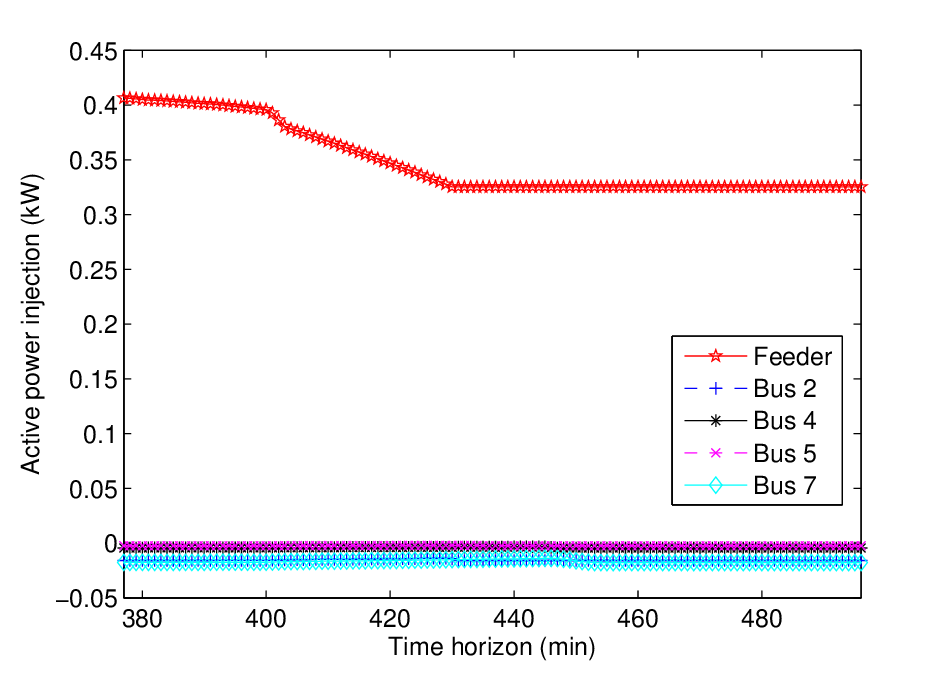}
\label{fig:activeplot_2hr}}
\subfigure[Reactive power injections.]{
\includegraphics[scale=0.25]{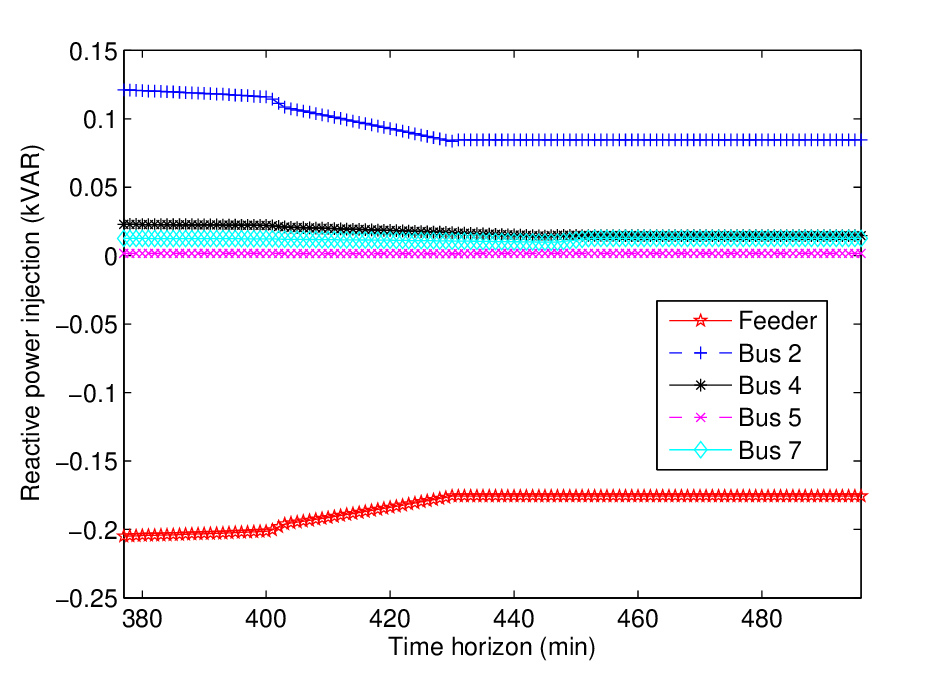}
\label{fig:reactiveplot_2hr}}
\caption{Active and reactive power injections at various buses in the 34-bus network.}
\label{fig:2hr}
\end{figure}
\begin{figure}[t!]
  \centering
  \subfigure[34-bus system]{
  	 \includegraphics[width=4.1cm]{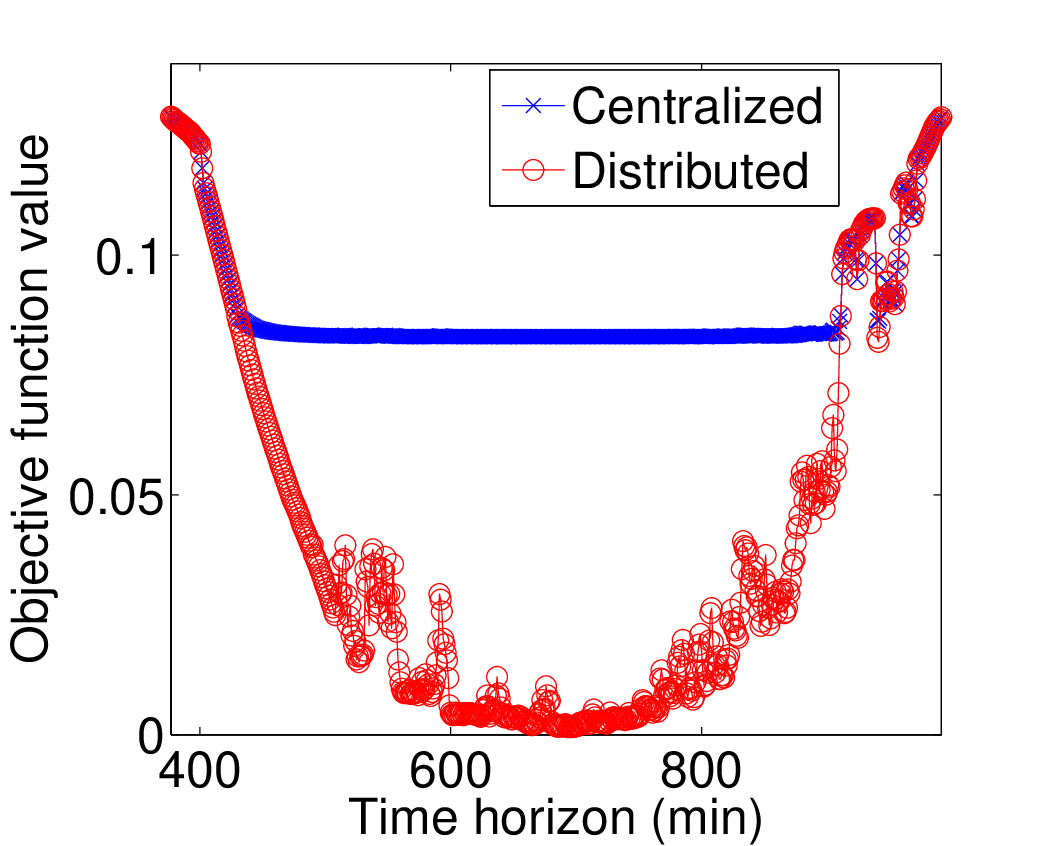}
	 }
  \subfigure[123-bus system]{
  	\includegraphics[width=4.1cm]{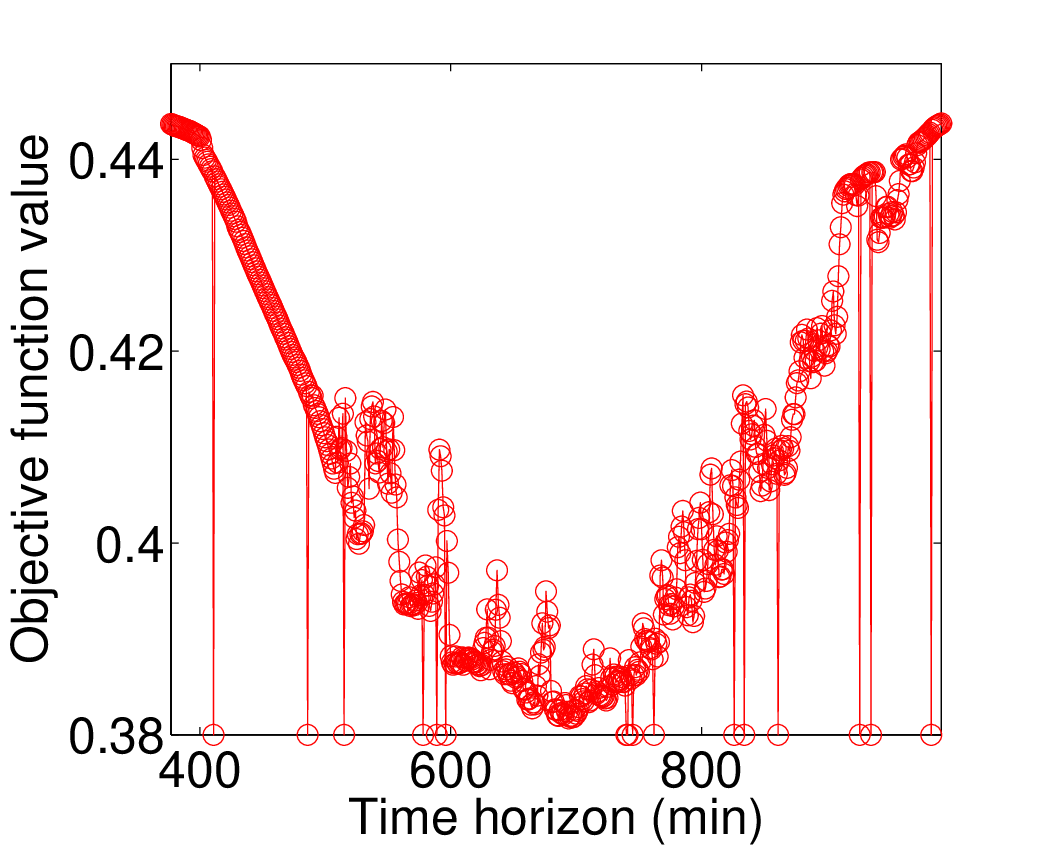} \label{fig:122_obj_plot}
	}
    \caption{Objective function values  computed by the distributed algorithm.}
    \label{fig:obj_plot}
\end{figure}

In order to check if the distributed algorithm can achieve the global optimum, we compare the objective function values computed by the distributed algorithm to those by the centralized solver (e.g., SDP3 or SeDuMi \cite{S98guide});  the results are plotted in Fig.~\ref{fig:obj_plot}. The active and reactive power injection at various buses in the 34-bus network are shown in Fig. \ref{fig:activeplot_2hr} and Fig. \ref{fig:reactiveplot_2hr} respectively. 
{As observed in Fig. \ref{fig:obj_plot},} for the 34-bus system, we can see that the distributed algorithm converges to the optimum all the time, whereas for the 123-bus system convergence occurs most of the time. {In Fig. \ref{fig:122_obj_plot}, the dropping lines correspond to the non-convergent cases where the convergence fails because the pre-defined 300 iterations allowed  were exhausted. In this case, the value from the previous solution is used. As shown in the simulation, the voltages are still maintained at their reference values.}
 In the 34-bus and 123-bus systems, the centralized solver failed to solve the system due to convergence issues. 
 
\begin{figure}[t!]
\centering
\subfigure[34-Bus Network]{
\includegraphics[scale=0.25]{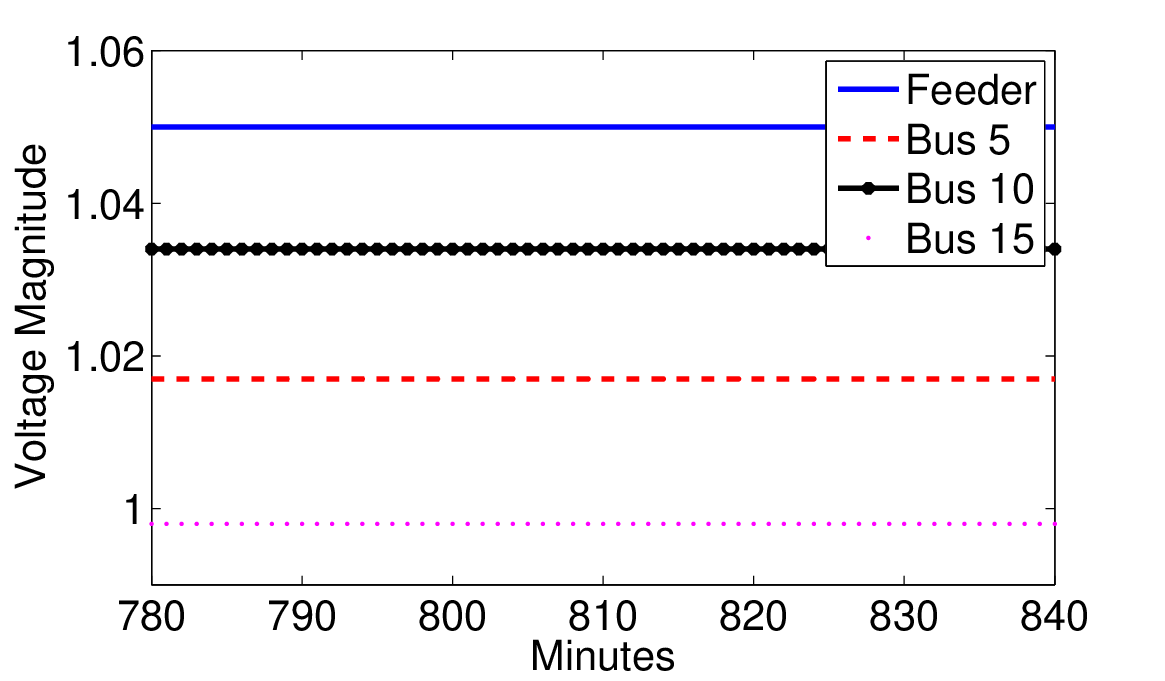}
\label{fig:voltages_34}}
\subfigure[123-Bus Network]{
\includegraphics[scale=0.25]{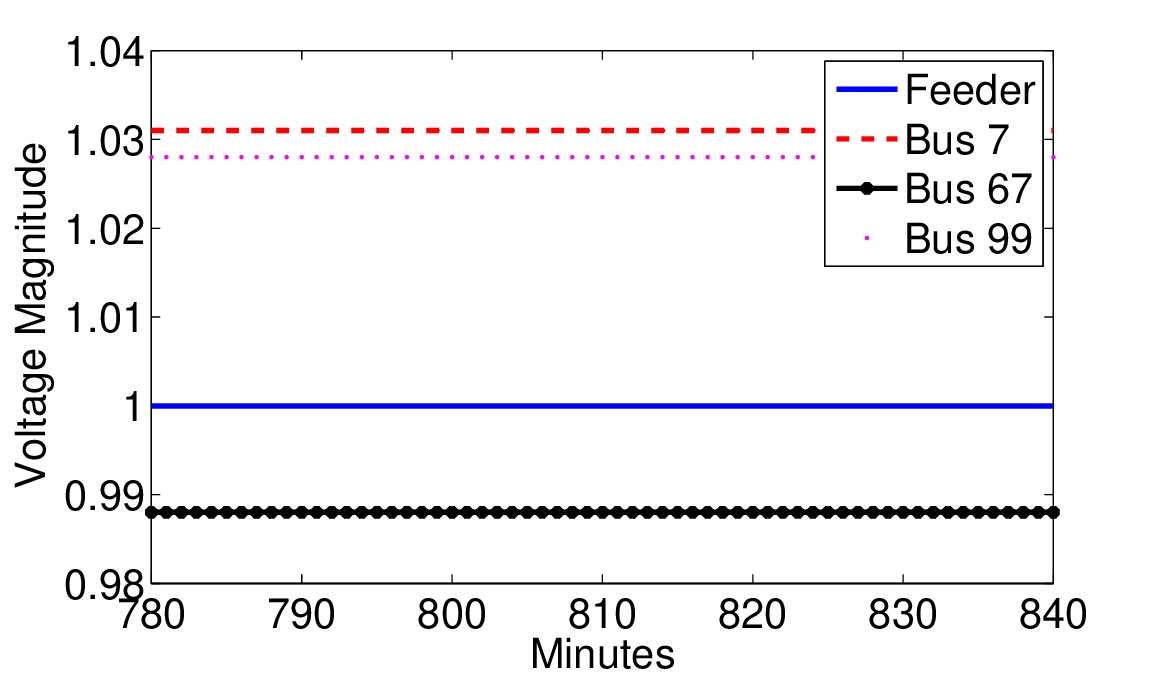}
\label{fig:voltages_123}
}
\caption{Voltage profile over time at representative buses. The proposed voltage regulation method is able to keep the voltages constant at their perspective references values. }
\label{fig:voltages}
\vspace{-0.15in}
\end{figure}

 {
Fig.~\ref{fig:voltages} displays the voltage profile at various representative buses of the the  34- and 123-bus test system over a one-hour period with high variability in the the $P^{PV}_i$'s caused by the  high-variability irradiance period displayed in Fig.~\ref{hour_horizon}.  This one-hour period  corresponds to the portion of the daily irradiance profile in  Fig~ \ref{entire_horizon} between the $781^{th}$ and $840^{th}$ minutes. For this simulation, the settings of the  conventional voltage regulation devices are kept at the values given in  \cite{testfeeders}, whereas the  $V^{ref}_i$'s  in \eqref{eqn:vol} result from the solution to the power flow equations for the nominal $\hat{P}_i$'s as specified in  \cite{testfeeders}. The fact that all the voltages displayed in Fig.~\ref{fig:voltages} remain at their reference value illustrates the effectiveness of our proposed voltage regulation method to mitigate the effect of fast-varying power injections arising from PV systems.  If a system do not use the reactive capability of the DERs and experiences high penetrations of solar-based generation, the voltages could exceed design tolerances by fluctuating outside the prescribed magnitude interval [0.95,1.05 p.u.].
} 
%
%
%
%
 
\begin{figure}[t!]
\vspace{0.1in}
  \centering
  \subfigure[34-bus system]{
    \includegraphics[width=4.1cm]{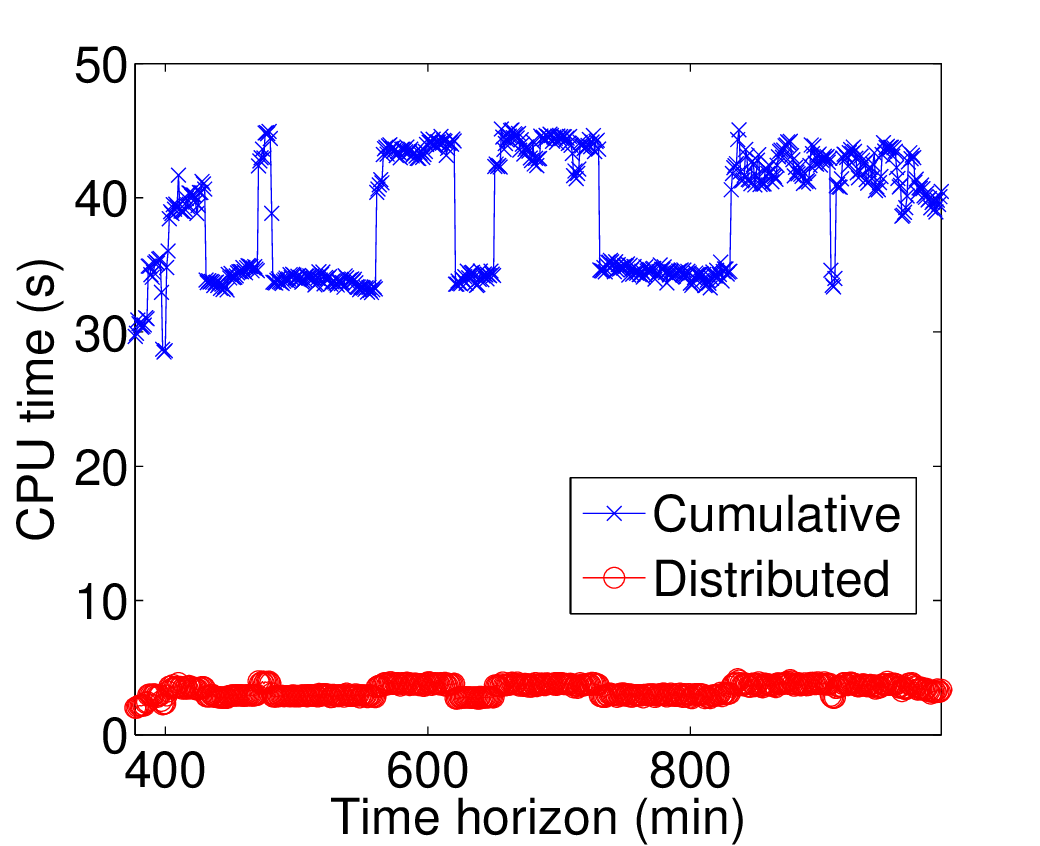}
    }
  \subfigure[123-bus system]{
  	\includegraphics[width=4.1cm]{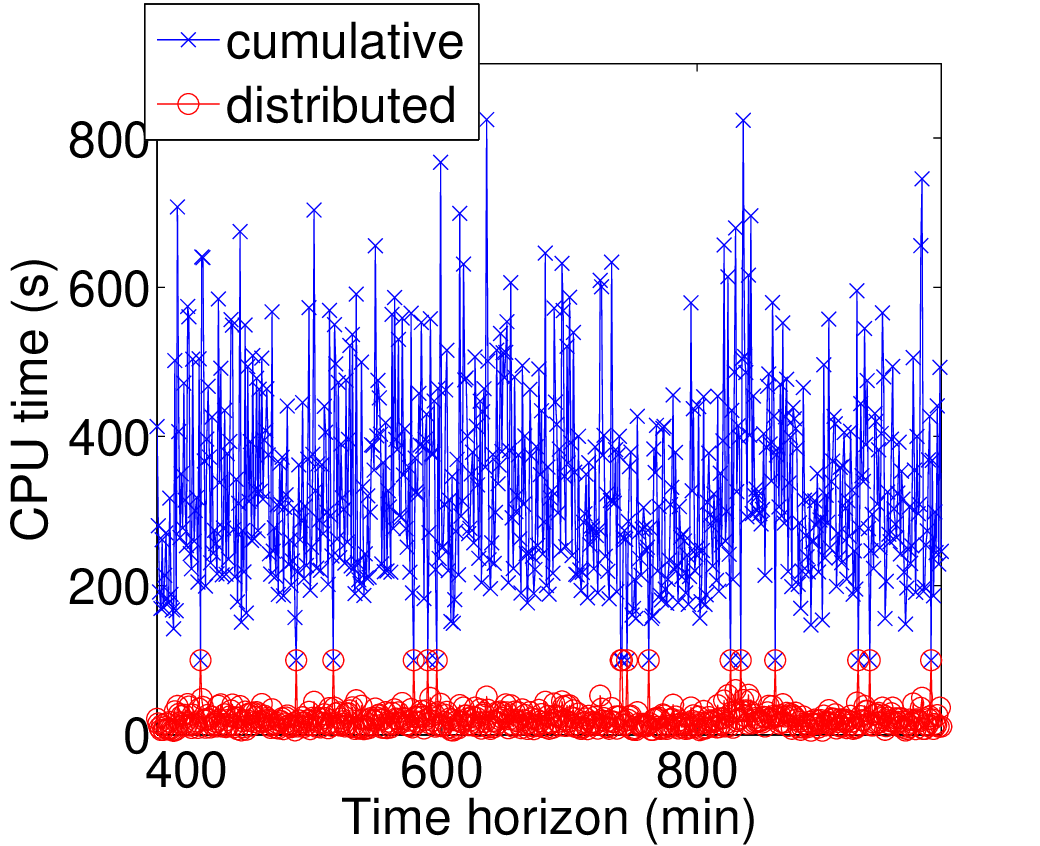}
	}
    \caption{Computation time of the distributed algorithm.}
    \label{fig:timeplot}
    \vspace{-0.45in}
\end{figure}

Fig.~\ref{fig:timeplot} shows the computational times corresponding to each test system; here we only consider the CPU time spent on the SDP solver and assume that  communication overheads can be neglected. In our simulation, we implement the algorithm iteratively; in each iteration, we solve the subproblems sequentially. In Fig.~\ref{fig:timeplot}, each subfigure contains two curves.  One (distributed) is to sum the CPU times of the subproblems which need the longest CPU time in each iteration. In other words, we only consider the most demanding subproblem in each iteration and then sum the CPU times spent on these subproblems in all iterations. 
The average CPU computation time  for the three cases are $1.28$ s, $3.33$ s, and $19.69$ s, respectively, which are substantially shorter than the one-minute cycles considered.

\begin{figure}[t!]
 \vspace{-0.15in}
\centering
\includegraphics[scale=0.39]{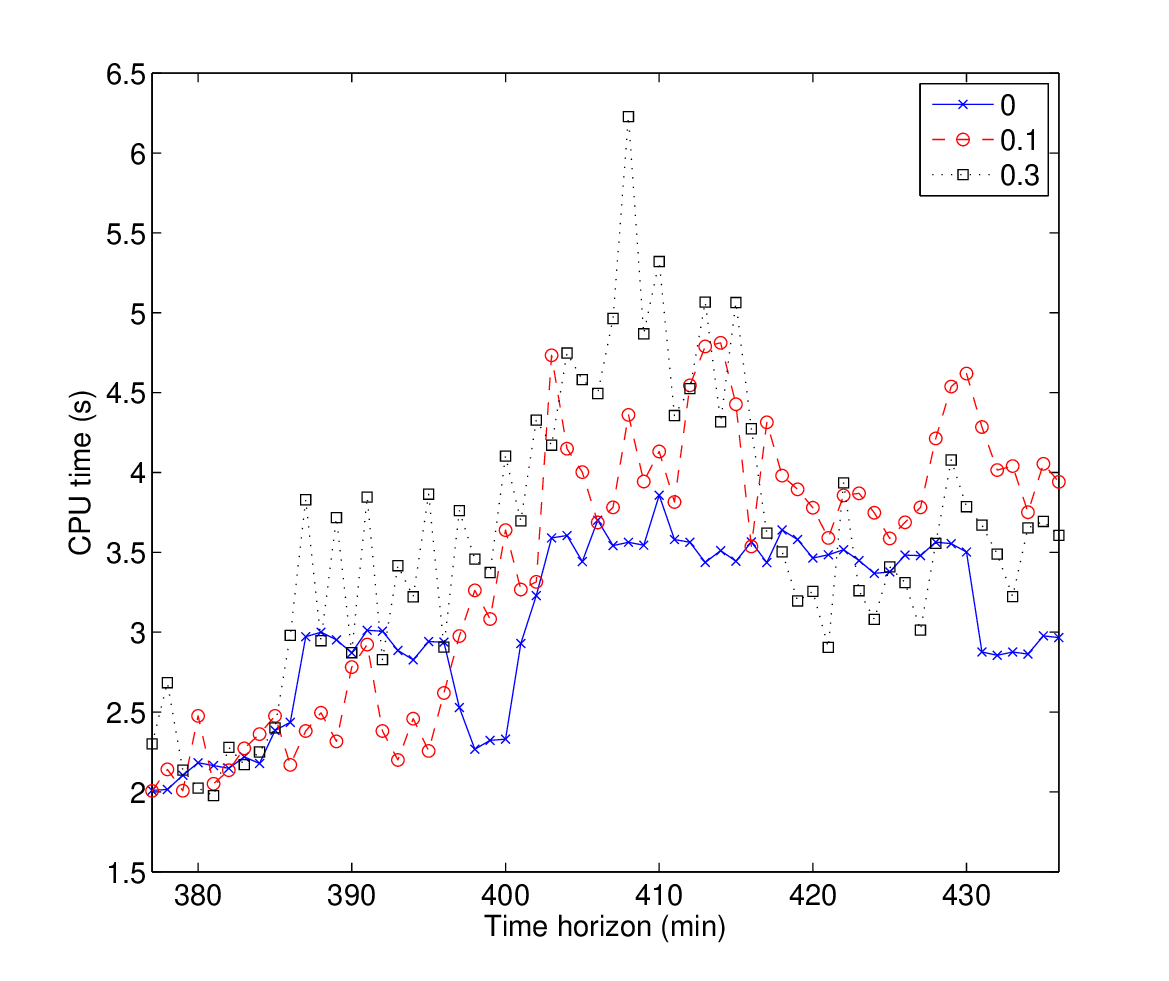}
\caption{Time it takes for Algorithm~\ref{algorithm_1} to converge under the presence of communication link failures.}
\label{fig:robust}
  \vspace{-0.1in}
\end{figure}

Next, we show that the distributed algorithm is robust against  random communication link failures. We model communication failures as packets drops. This means that, at a given iteration, the Lagrangian multiplier transmitted on any particular edge could be lost with probability $p$, independent of all other transmissions.  Figure \ref{fig:robust} shows the average time of convergence needed over the day for the 34-bus network for $p=0$, $p=0.1$ and $p=0.3$;  convergence is always achieved.

\textbf{Remark~1.}
The  results displayed in Fig.~\ref{fig:obj_plot} correspond to a centralized solution of  \eqref{eqn:W}; however, we also performed simulations using a centralized solver (SeDuMi) to obtain a solution to \eqref{prob3}. However, the algorithm failed to converge most of the time even for the 34-bus network. We suspect that since there are duplicated variables in \eqref{prob3}, a naive implementation would be rather inefficient; however, a more careful centralized implementation (using, e.g., an SOCP formulation)   is likely speed up   \eqref{prob3}.  \qed

\section{Concluding Remarks} \label{sec:con}
We proposed a convex optimization based method to solve the voltage regulation problem in distribution networks. We cast the problem as a loss minimization program.  We showed that under broad conditions that are likely to be satisfied in practice, the optimization problem can be solved via its convex relaxation. We then proposed a {distributed algorithm that can be implemented in a network with a large number of buses; we demonstrated the effectiveness and robustness of the algorithm with two case studies.}

{As noted earlier, the proposed voltage regulation method is intended  to supplement the action of conventional voltage regulation devices. In this regard, throughout the paper, we assumed that there is a separation in the (slow) time-scale  in which the settings of   conventional voltage regulation devices are adjusted and the (fast) time-scale in which our proposed method operates.  With respect to this, a research direction worthy exploring is to carefully consider the coupling across the two time-scales and study the interplay between the optimal use of conventional voltage regulation devices in longer time-scales, and the use of our voltage regulation  method in shorter time-scales. This would allow us to, e.g., study the trade-offs between the location and number of conventional voltage regulation devices, and the location and number of DERs and DRRs with capability of providing reactive power. }

\appendix
\begin{IEEEproof}[Proof of Theorem \ref{thm:main}]
The first and third cases are clear. The interesting case to prove is to show that if the matrix $\bd{W}^*$ has rank higher than $1$, there is no rank-1 matrix $\bd{W}$ that results in a  feasible solution.

The requirement that $\ul{Q}_i < \be_i$ for $i=2,\dots,n$ is to ensure that the reactive lower bound is in fact never tight for all the nodes in the network. Let $h$ be the parent of $i$ and $k$ be a child of $i$. Since we assume that power always flow from parents to children in the network, and from the angle constraint in \eqref{eqn:ang_con}, $0 \leq \te_{hi} \leq \tan^{-1}(\frac{b_{ik}}{g_{ik}})$. Over this range, $Q_{ih} \geq 0$. This corresponds to the intuition that reactive power should flow up the tree to support the voltage. Note the inductive line is very lossy in terms of reactive powers, therefore $i$ might receive or supply reactive power to $k$. The $Q_{ik}$ is monotonic in $\te_{ik}$ starting at $\te_{ik}=0$ until it reaches is minimum at an angle of $ \tan^{-1} (\frac{g_{ik}}{b_{ik}})$. Let $\tilde{\te}_{ik}=\min(\tan^{-1} (\frac{g_{ik}}{b_{ik}}), \ov{\te}_{ik})$, then
$
Q_i  = Q_{ih} + \sum_{k: k \in \mc{C}(i)} Q_{ik}  \geq  \sum_{k: k \in \mc{C}(i)} Q_{ik}   \geq \sum_{k: k \in \mc{C}(i)} b_{ik}-g_{ik} \sin (\tl{\te}_{ik})- b_{ik} \cos (\te(\te)_{ik})  = \be_i.
$
Thus, if $\ul{Q}_i < \be_i$, the lower bounds on the reactive power injections are never tight.

To finish the proof we need to introduce some notations from \cite{LTZ12}.
For a $n$ bus network, let $\{\ov{\te}_{ik}\}$ be the set of angle constraints, one for each line. Then, the angle-constrained {\em active power injection region} is the set of all active  power injection vectors that satisfy the line angle constraints, i.e.,
$\mc{P}_\te = \{ \bd{p}  : \bd{p}=\Real \{\diag(\bd{v}\bd{v}^H \bd{Y}^H) \},~ |V_i|=1, |\te_{ik}| \leq \ov{\te}_{ik}\}$.
Let $\mc{F}_\te \subset \R^{2n-2}$ be the Cartesian product region of the $n-1$ active line flow regions, then
$
\mc{F}_\te= \Pi_{i \sim k} \mc{F}_{\te, ik}$.
Let $\bd{M} \in \R^{n \times 2n-2}$ be a matrix with the rows indexed by the buses and the columns indexed by the $2(n-1)$ ordered pair of edges, i.e., if $i$ is connected to $k$, both $(i,k)$ and $(k,i)$ are included; thus $\bd{M}[i,(k,l)]=1$ if $i=k$ or $i=l$, and $\bd{M}[i,(k,l)]=0$ otherwise.
$\bd{M}$ is a generalized edge to bus incidence matrix, and $\mc{P}_\te= \bd{M} \mc{F}_\te$ i.e., the power injection region is obtained by a linear transformation of the product of line flow regions.

Similarly, $\mc{G}_\te$ is the product region of the $n-1$ reactive line flow regions. Then, for all $(i,k)$, by stacking the   $\bd{H}_{ik}$'s as defined in \eqref{transformation} into a $2(n-1) \times 2(n-1)$ block diagonal matrix $\bd{H}$, i.e., $\bd{H}= \diag(\{H_{ik}\}_{i \sim k})$, we obtain a   the global transform between $\mc{F}$ and $\mc{G}$. The angle-constrained {\em reactive power injection region} $\mc{Q}_\te$ is given by $\mc{Q}_\te = \bd{M} \mc{G}_\te= \bd{M} \bd{H} \mc{F}_\te$.

Since by construction the lower bonds on the reactive power injection are never tight, we can  ignore them from now on. Let $\mc{P}$ be the feasible region of the original problem \eqref{eqn:v}, that is,
$\mc{P}=\{ \bd{p}: \exists \bd{v} \in \C^n, P_i = \Tr (\bd{A}_i \bd{v}\bd{v}^H), |V_i|=1, \ul{P}_i \leq P_i \leq \ov{P}_i, \Tr(\bd{B}_i \bd{v}\bd{v}^H) \leq \ov{Q}_i, |\te_{ik}| < \ov{\te}_{ik}, \forall i \sim k\}$.
We can equivalently write $\mc{P}$ as
$
\mc{P}=\bd{M} (\mc{F}_\te \cap \mc{F}_P \cap \mc{F}_Q),
$
where  $\mc{F}_P$ is the flow region satisfying the real power constraints, that is, $\mc{F}_P = \{ \bd{f} \in \R^{2n-2}: \bd{p}=\bd{M} \bd{f}, \ul{P}_i \leq P_i \leq \ov{P}_i \}$. $\mc{F}_Q$ is the flow region satisfying the reactive power constraints, that is, $\mc{P}_Q=\{\bd{f} \in \R^{2n-2}: \bd{q}=\bd{M} \bd{H}\bd{f}, Q_i \leq \ov{Q}_i \}$. Since $\mc{F}_P$ and $\mc{F}_Q$ are defined by linear inequalities, they are convex. However, $\mc{F}_\te$ is not.

Let $\mc{S}$ be the feasible region of the relaxed problem \eqref{eqn:W}.
It turns out that $\mc{S}= \bd{M} (\convhull(\mc{F}_\te) \cap \mc{F}_P \cap \mc{F}_Q)$, is convex, and contains $\mc{P}$.

Now we need to define the Pareto-front of a set. Let $\mc{X} \subset \R^n$, we say $\x \in \mc{X}$ is Pareto-optimal if $\not\exists \bd{y} \in \mc{X}$ such that $\bd{y} \leq \bd{x}$ with strict inequality in at least one coordinate. The set of Pareto-optimal points is called the Pareto-front of $\mc{X}$, and labeled $\mc{O}(\mc{X})$. When minimizing a strictly increasing function, the optimal is always achieved in the Pareto-front. Therefore to show the second statement in the theorem, it suffices to show the following lemma.
\begin{lem}
Suppose $\mc{P}$ is not empty, then $\mc{P}=\mc{O}(\mc{S})$.
\end{lem}
Suppose the lemma is true, then if the optimal solution of the relaxed problem \eqref{eqn:W} is of rank 2, then $\mc{P}$ must be empty.

The proof of this lemma is similar to Lemma~4 in \cite{LTZ12}. Let $\bd{p}^* \in \mc{S}$ be the optimal solution of the relaxed problem, $\bd{f}^* \in  \convhull(\mc{F}_\te) \cap \mc{F}_P \cap \mc{F}_Q$ its corresponding active flow vector and $\bd{r}^* = \bd{H} \bd{f}^* $ be the corresponding reactive power flow vector. It suffices to show that if $P_i^* > \ul{P}_i$, then $(f_{ik}^*,f_{ki}^*) \in \mc{F}_{\te,ik}$ for every $k \sim i$. Once this fact is established, the rest of the proof is the same as the proof of Lemma 4 in \cite{LTZ12}. Suppose that $P_i^* > \ul{P}_i$, but $(f_{ik}^*,f_{ki}^*) \notin \mc{F}_{\te,ik}$ for some $k$. Then there exists $\epsilon >0$ such that $(f_{ik}^*-\epsilon,f_{ki}^*) \in \convhull(\mc{F}_{\te,ik})$. Let $(\tl{f}_{ik},\tl{f}_{ki})= (f_{ik}^*-\epsilon,f_{ki}^*)$. Since
\begin{equation*}
\bd{H}_{ik} \bma -\epsilon \\ 0 \ebma =
-\frac{\epsilon}{2 b_{ik} g_{ik}} \bma b_{ik}^2-g_{ik}^2 \\ b_{ik}^2 + g_{ik}^2 \ebma < \bd{r}^*,
\end{equation*}
Therefore $(\tl{f}_{ik},\tl{f}_{ki})$ is a better feasible flow on the line $(i,k)$, which contradicts the optimality of $\bd{f}^*$.
\end{IEEEproof}

\vspace{-0.1in}
\bibliographystyle{IEEEtran}
\bibliography{IEEEabrv,mybibv2}

\begin{IEEEbiography}[{\includegraphics[width=1in,height=1.25in,clip,keepaspectratio]{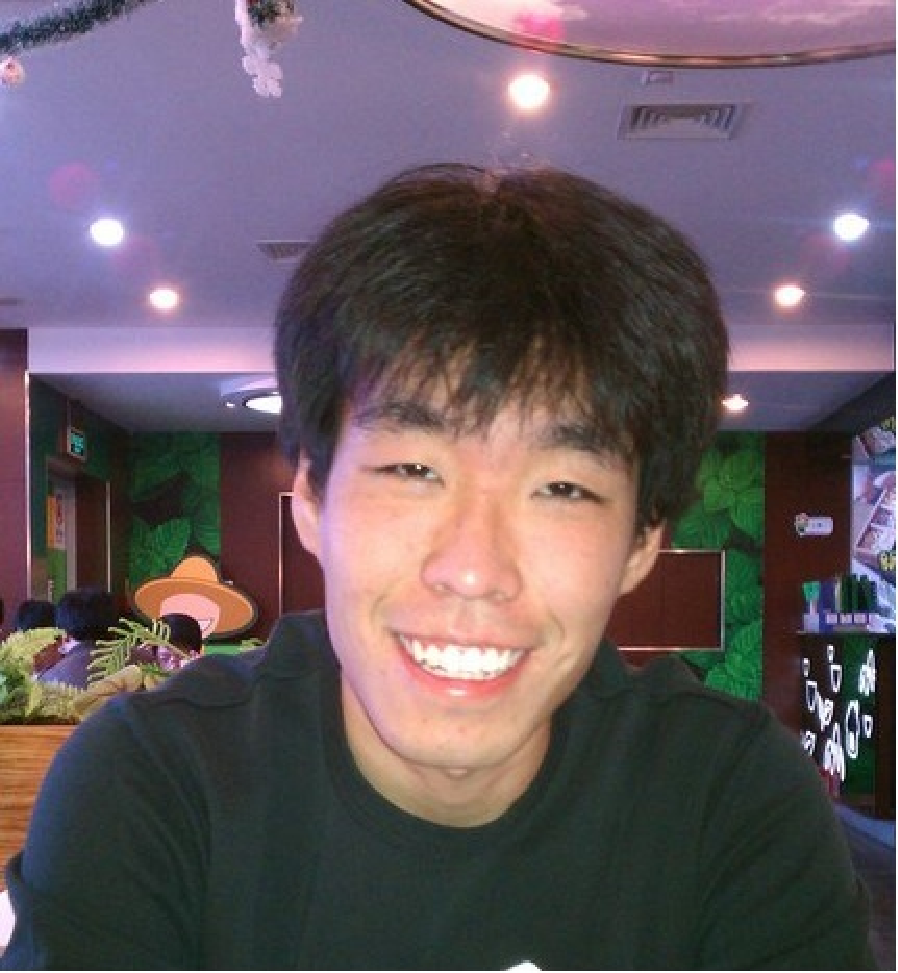}}]{Baosen Zhang} is a postdoctoral scholar at Stanford University, jointly hosted by departments of Civil and Environmental Engineering and Management \& Science Engineering. He will start as an Assistant Professor in Electrical Engineering at the University of Washington in 2015. 

He received the B.A.Sc. degree in engineering science from the University of Toronto, Toronto, ON, Canada, in 2008  and the Ph.D. degree Department of Electrical Engineering and Computer Science, University of California at Berkeley in 2013. His interest is in the area of power systems, particularly in the  fundamentals of power flow and the economical challenges resulting from renewables. He was awarded several fellowships: the Post Graduate Scholarship from NSERC in 2011; the Canadian Graduate Scholarship from NSERC in 2008; the EECS fellowship from Berkeley in 2008.  
\end{IEEEbiography}

\begin{IEEEbiography}[{\includegraphics[width=1in,height=1.25in,clip,keepaspectratio]{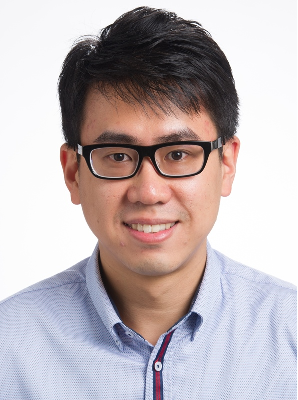}}]{Albert Y.S. Lam}
 received the BEng degree (First Class Honors) in Information Engineering and the PhD degree in Electrical and Electronic Engineering from the University of Hong Kong, Hong Kong, China, in 2005 and 2010, respectively. He is a research assistant professor at the Department of Computer Science of Hong Kong Baptist University, Hong Kong and he was a postdoctoral scholar at the Department of Electrical Engineering and Computer Sciences of University of California, Berkeley, CA, USA, in 2010--12. He is a Croucher research fellow. His research interests include optimization theory and algorithms, evolutionary computation, smart grid, smart city planning, and Internet protocols and applications.
\end{IEEEbiography}

\begin{IEEEbiography}[{\includegraphics[width=1in,height=1.25in,clip,keepaspectratio]{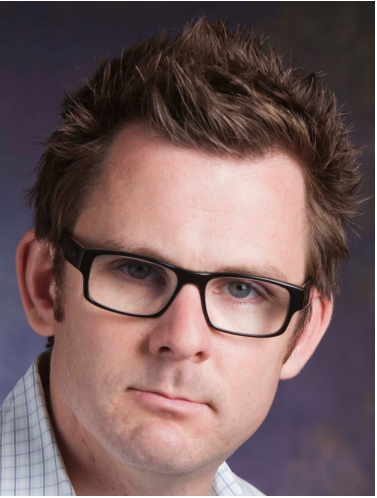}}]{Alejandro D. Dom\'{i}nguez-Garc\'{i}a (S'02--M'07)} is an Assistant Professor in the Department of Electrical and Computer Engineering of the University of Illinois at Urbana-Champaign, where he is affiliated with the Power and Energy Systems area. His research interests are in the areas of system reliability theory and control, and their application to electric power systems, power electronics, and embedded electronic systems for safety-critical/fault-tolerant aircraft, aerospace, and automotive applications. He received the Ph.D. degree in Electrical Engineering and Computer Science from the Massachusetts Institute of Technology, Cambridge, MA, in 2007 and the degree of Electrical Engineer from the University of Oviedo (Spain) in 2001. After finishing his Ph.D., he spent some time as a post-doctoral research associate at the Laboratory for Electromagnetic and Electronic Systems of the Massachusetts Institute of Technology. Dr. Dom\'{i}nguez-Garc\'{i}a received the NSF CAREER Award in 2010, and the Young Engineer Award from the  IEEE Power and Energy Society in 2012. He is also a Grainger Associate since 2011. He currently serves as an Associate Editor for the IEEE Transactions on Power Systems and the IEEE Power Engineering Letters.
\end{IEEEbiography}

\begin{IEEEbiography}[{\includegraphics[width=1in,height=1.25in,clip,keepaspectratio]{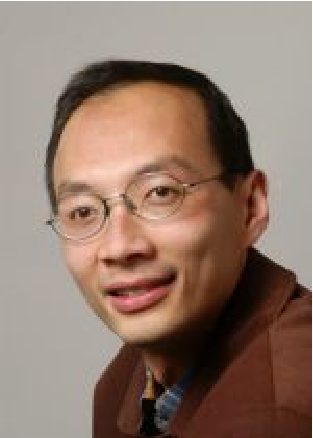}}]{David Tse} received the B.A.Sc. degree in systems design engineering from University of Waterloo in 1989, and the M.S. and Ph.D. degrees in electrical engineering from Massachusetts Institute of Technology in 1991 and 1994 respectively. From 1994 to 1995, he was a postdoctoral member of technical staff at A.T. \& T. Bell Laboratories. From 1995- 2014, he was on the faculty of the University of California at Berkeley. He is currently a professor at Stanford University.  He received a 1967 NSERC graduate fellowship from the government of Canada in 1989, a NSF CAREER award in 1998, the Best Paper Awards at the Infocom 1998 and Infocom 2001 conferences, the Erlang Prize in 2000 from the INFORMS Applied Probability Society, the IEEE Communications and Information Theory Society Joint Paper Awards in 2001 and 2013, the Information Theory Society Paper Award in 2003, the 2009 Frederick Emmons Terman Award from the American Society for Engineering Education, a Gilbreth Lectureship from the National Academy of Engineering in 2012, the Signal Processing Society Best Paper Award in 2012 and the Stephen O. Rice Paper Award in 2013. He is a coauthor, with Pramod Viswanath, of the text "Fundamentals of Wireless Communication", which has been used in over 60 institutions around the world.
\end{IEEEbiography}

\end{document}